\documentclass[12pt,leqno]{amsart}

\usepackage[T1]{fontenc}
\usepackage[utf8]{inputenc}
\usepackage[english]{babel}

\usepackage{amsmath,amsthm,amsfonts,amssymb}
\usepackage{dsfont}
\usepackage{esint}
\usepackage{mathrsfs}
\usepackage{cancel}
\usepackage{float}

\usepackage{verbatim}
\usepackage{xurl}
\usepackage[numbers,sort&compress]{natbib}
\usepackage{needspace}
\usepackage[colorlinks=true]{hyperref}
\usepackage[dvipsnames]{xcolor}
\hypersetup{
    colorlinks=true,
    linkcolor=PineGreen,
    citecolor=PineGreen,
    urlcolor=NavyBlue
}
\DeclareUrlCommand\doiurl{\urlstyle{rm}}

\usepackage[a4paper,margin=1in]{geometry}
\usepackage{graphicx}
\usepackage{pgfplots}
\usepgfplotslibrary{groupplots}
\pgfplotsset{compat=1.16}

\pgfmathdeclarefunction{hbump}{1}{\pgfmathparse{#1>0.02 ? exp(-1/(#1)) : 0}}
\pgfmathdeclarefunction{dhbump}{1}{\pgfmathparse{#1>0.02 ? exp(-1/(#1))/((#1)*(#1)) : 0}}
\pgfmathdeclarefunction{chicut}{1}{%
  \pgfmathparse{hbump(1-abs(#1))/(hbump(1-abs(#1))+hbump(abs(#1)-0.5))}}
\pgfmathdeclarefunction{dchicut}{1}{%
  \pgfmathparse{-sign(#1)*(dhbump(1-abs(#1))*hbump(abs(#1)-0.5)+hbump(1-abs(#1))*dhbump(abs(#1)-0.5))/(hbump(1-abs(#1))+hbump(abs(#1)-0.5))^2}}

\pgfmathsetmacro{\plotM}{400}
\pgfmathsetmacro{\plotMs}{sqrt(\plotM)}
\pgfmathsetmacro{\plotX}{1/\plotMs}

\pgfmathdeclarefunction{hsafe}{1}{\pgfmathparse{(#1>0.02)*exp(-1/max(#1,0.02))}}
\pgfmathdeclarefunction{chisafe}{1}{%
  \pgfmathparse{hsafe(1-abs(#1))/(hsafe(1-abs(#1))+hsafe(abs(#1)-0.5))}}

\pgfmathdeclarefunction{bgw}{1}{\pgfmathparse{-0.75*chisafe(#1)*sin(deg(3*pi*(#1)))}}

\pgfmathdeclarefunction{ptw}{3}{\pgfmathparse{-(#2)*chisafe((#1)/(#3))*sin(deg(4*pi*(#1)/(#3)))}}

\theoremstyle{plain}
\newtheorem{thm}{Theorem}[section]
\newtheorem{cor}[thm]{Corollary}
\newtheorem{lem}[thm]{Lemma}

\theoremstyle{definition}
\newtheorem{dfn}[thm]{Definition}

\theoremstyle{remark}
\newtheorem{rem}[thm]{Remark}

\newcommand{\R}{\mathbb{R}}

\newcommand{\N}{\mathbb{N}}
\newcommand{\T}{\mathbb{T}}

\newcommand{\dd}{\text{d}}

\newcommand{\supp}{\operatorname{supp}}

\newcommand{\ind}{\mathds{1}}

\begin{document}

\title[A gentle introduction to the vortex layer cascade method]{A gentle introduction to the vortex layer cascade method: finite-time blow-up for the forced De Gregorio equation}
\author{Diego Córdoba}
\address{Instituto de Ciencias Matemáticas CSIC-UAM-UCM-UC3M, Madrid, Spain}
\email{dcg@icmat.es}

\author{José Antonio Fuentes-Mesa}
\address{Instituto de Ciencias Matemáticas CSIC-UAM-UCM-UC3M, Madrid, Spain}
\email{jose.fuentes@icmat.es}

\author{Luis Martínez-Zoroa}
\address{CUNEF Universidad, Madrid, Spain}
\email{luis.martinezzoroa@cunef.edu}
\date{}

\begin{abstract}
This paper offers a detailed introduction to the vortex layer cascade method through a toy model, the forced De Gregorio equation on the real line. We show that for every $\delta\in(0,1)$, there exist a forcing term $F$, supported in $[-1,1]$ and bounded in $C^{0,1-\delta}$ uniformly up to and including the blow-up time $t=1$, and a corresponding solution $w$, smooth and supported in $[-1,1]$ for $t\in[0,1)$, such that $\limsup_{t\nearrow1}\|w(\cdot,t)\|_{L^\infty}=\infty$. The singularity is produced by an infinite cascade of highly oscillatory layers, each of which is stretched and amplified by the flow of the previous ones.
\end{abstract}

\maketitle

\tableofcontents

\section{Introduction}
\subsection{Context}
Recent developments concerning finite-time singularity formation in fluid equations have attracted considerable attention, particularly following three works by Alpöge, Buckmaster, and collaborators establishing finite-time blow-up with smooth forcing for three distinct fluid models: the three-dimensional incompressible Euler equations~\cite{AlpogeBuckmaster2026Euler}, the inviscid Boussinesq system~\cite{AlpogeBuckmaster2026Boussinesq}, and the incompressible porous media equation on the torus~\cite{AlpogeBuckmasterCoiculescu2026}. All three works employ the vortex layer cascade method introduced by Córdoba and Martínez-Zoroa~\cite{CordobaMartinezZoroaBlowupEuler} and subsequently developed in related works by Córdoba, Martínez-Zoroa, and collaborators~\cite{Andres1, Andres2, CordobaMartinezZoroaZhengBlowup, JoseAntonio1, CordobaMartinezZoroaBlowupIPM}. See also the seminar notes~\cite{lain2026notes} for a slightly different description of the method and a sketch of its application to the two-dimensional Boussinesq and non-homogeneous Euler equations.

Motivated by these results and preliminary work developed by one of the authors in their master's thesis~\cite{TFM}, in this paper we study the forced De Gregorio equation on the real line,
\begin{equation}\label{eq:DGF}\tag{DGF}
    \begin{cases}
        \partial_t w+u(w)\partial_xw=wHw+F, \quad \text{in }\R\times(0,\infty), \\
        u(w)(x,t)=\int_0^xHw(y,t)\dd y, \quad  \text{in }\R\times[0,\infty),\\
        w(x,0)=w_0(x), \quad \text{in }\R,
    \end{cases}
\end{equation}
as a one-dimensional setting in which to investigate the vortex layer cascade method. Here $H$ denotes the Hilbert transform on $\R$ (see Subsection~\ref{sec:prelim} for our normalization) and $F$ is a forcing term. The De Gregorio equation originates in the study of vortex stretching for the incompressible three-dimensional Euler equations, where the vorticity $\omega$ satisfies
$\partial_t\omega+u\cdot\nabla\omega=\omega\cdot\nabla u$.

In 1985, Constantin, Lax and Majda~\cite{ConstantinLaxMajda1985} introduced a one-dimensional model designed to capture the essential nonlocal structure of the vortex-stretching mechanism:
\begin{equation}\label{eq:CLM}\tag{CLM}
    \begin{cases}
        \partial_t w = wHw, \quad  \text{in }\R\times(0,\infty),\\
        w(x,0)=w_0(x), \quad \text{in }\R.
    \end{cases}
\end{equation}
This model exhibits finite-time blow-up for a suitable class of nontrivial initial data. 

In 1990, De Gregorio~\cite{DeGregorio1990} proposed a modification incorporating an advective term, thereby accounting for the interplay between transport and vortex stretching:
\begin{equation}\label{eq:DG}\tag{DG}
    \begin{cases}
    \partial_t w+u\partial_xw=wHw, \quad \text{in }\R\times(0,\infty), \\
        \partial_xu=Hw, \quad  \text{in }\R\times[0,\infty),\\
        w(x,0)=w_0(x), \quad \text{in }\R,
    \end{cases}
\end{equation}
which is precisely Equation~\eqref{eq:DGF} with $F=0$ up to the normalization $u(0,t)=0$ determining the velocity $u$. Note that the relation $\partial_xu=Hw$ mimics the Biot--Savart law relating the velocity and the vorticity in the incompressible three-dimensional Euler equations.

Okamoto, Sakajo and Wunsch~\cite{OkamotoSakajoWunsch2008} subsequently considered the one-parameter family
\begin{equation}\label{eq:gCLM}\tag{gCLM}
\begin{cases}
    \partial_t w+au\partial_xw=wHw, \quad \text{in }\T\times(0,\infty), \\
        \partial_xu=Hw, \quad  \text{in }\T\times[0,\infty),\\
        w(x,0)=w_0(x), \quad \text{in }\T,
    \end{cases}
\end{equation}
for each $a\in\R$, which interpolates between the Constantin--Lax--Majda model ($a=0$) and the De Gregorio model ($a=1$) in the periodic setting. For this family they established local existence of solutions and a Beale--Kato--Majda-type continuation criterion, and their numerical experiments suggested that sufficiently strong advection could prevent the formation of singularities.

Since then, the central question has been whether the competition between stretching and advection leads to finite-time blow-up.

When $a<0$, the transport and vortex stretching terms cooperate to promote singularity formation. Indeed, Castro and Córdoba~\cite{CastroCordoba2010} proved finite-time blow-up for smooth initial data in this regime.

For $a>0$, the question is more delicate because advection can have a stabilizing effect. Nevertheless, Elgindi and Jeong~\cite{ElgindiJeong2020} showed that this effect does not preclude finite-time singularity formation when $|a|$ is sufficiently small: they constructed smooth initial data leading to an exactly self-similar blow-up. Moreover, they proved that, for every $a\in\mathbb{R}$, there exists a sufficiently small Hölder exponent $\alpha>0$, depending on $a$, for which some initial data in $C^{0,\alpha}(\R)$ develop an exactly self-similar singularity in finite time.

Subsequently, Chen, Hou, and Huang~\cite{ChenHouHuang2021} established the existence of initial data in $C_c^\infty(\R)$ that lead to an expanding and asymptotically self-similar singularity in finite time for the original De Gregorio equation. Here, \emph{expanding} means that the scaling exponent is negative, which in particular implies that, although both the initial data and the self-similar profile are compactly supported, the support of the solution becomes unbounded as the blow-up time is approached. Moreover, the quotient $w(x,t)/x$ remains uniformly bounded up to the blow-up time, so the solution cannot blow up at any finite point, that is, the singularity forms at spatial infinity. They also proved finite-time blow-up for the generalized De Gregorio model for all $a\in\R$ on the torus, starting from compactly supported Hölder-continuous initial data whose regularity exponent is sufficiently small.

Later, Huang, Tong, and Wei~\cite{HuangTongWei2023} showed that the De Gregorio equation on the real line admits infinitely many compactly supported, self-similar solutions that are distinct under rescaling and blow up in finite time. Their basic family consists of eigenfunctions of a compact self-adjoint operator, and the leading eigenfunction coincides, up to rescaling, with the self-similar profile obtained in~\cite{ChenHouHuang2021}. These solutions are also of expanding type.

More recently, Zheng~\cite{Zheng2023} established the existence of exactly self-similar solutions that blow up in finite time in $C^{0,\alpha}(\R)$ for every $\alpha\in(0,1)$, provided that $|a\alpha|$ is sufficiently small, thereby removing the restriction on the Hölder exponent in the earlier construction of Elgindi and Jeong.

The situation is different on the torus. Chen~\cite{Chen2025} studied the De Gregorio equation on $\mathbb{S}^1$ for $\pi$-periodic initial data in the class $X$ of odd functions with constant sign on $[0,\pi/2]$. The smooth initial data for which finite-time blow-up is established in~\cite{ChenHouHuang2021} share these sign and symmetry properties. Chen also established a one-point blow-up criterion and proved global well-posedness for initial data $w_0\in H^1\cap X$ with $w_0(x)x^{-1}\in L^\infty$, under an additional technical assumption. In particular, finite-time blow-up from smooth initial data in $X$ is ruled out on the torus. On the other hand, for every $0<\alpha<1$ and $s<3/2$ he constructed initial data in $X\cap C^{0,\alpha}\cap H^s\cap C^\infty(\mathbb{S}^1\setminus\{0\})$ whose solutions develop a singularity in finite time.

A different approach to singularity formation was introduced by Córdoba, Martínez-Zoroa, and Zheng~\cite{CordobaMartinezZoroaZheng2025}, who constructed finite-time singularities for both the generalized De Gregorio model and the three-dimensional incompressible Euler equations using infinitely many vorticity regions separated by vortex-free regions. In their construction, the entire hierarchy of vortex layers is encoded in the initial data, which makes this approach related to, but fundamentally different from the strategy developed by Córdoba and Martínez-Zoroa~\cite{CordobaMartinezZoroaBlowupEuler}.

In this paper, we investigate whether the vortex layer cascade method can be implemented for Equation~\eqref{eq:DGF}. We emphasize that finite-time blow-up for the unforced De Gregorio equation with smooth, compactly supported initial data is already known~\cite{ChenHouHuang2021}. The purpose of this paper is different: the singularity constructed here is produced by the vortex layer cascade, a non-self-similar mechanism of a different nature, and the one-dimensional setting allows us to present every step of the method explicitly.  

Indeed, the solution given by Theorem~\ref{construccionblowup} is qualitatively different from the expanding singularities of~\cite{ChenHouHuang2021,HuangTongWei2023}, as it remains supported in $[-1,1]$ for all times before the blow-up time (and so does the forcing term), so the singularity forms within $[-1,1]$. 

Although our solution does not have a constant sign on $(0,\infty)$ since the layers are highly oscillatory, the way in which the singularity forms is closely related to the one-point blow-up criterion of Chen~\cite{Chen2025}, as the stretching rate at the origin, $\partial_xu(w)(0,t)=Hw(0,t)$, grows very fast as the blow-up time is approached (and so does its time integral). 

\Needspace{9\baselineskip}
\subsection{Main result}
The main result of this paper is the following theorem.
\begin{thm}\label{construccionblowup}
For each $\delta\in(0,1)$, there exist a classical solution $w$ to Equation~\eqref{eq:DGF} and a forcing term $F$ such that
$$
w,F\in C^\infty(\mathbb{R}\times[0,1)),
\qquad
F\in C\bigl([0,1];C^{0,1-\delta}(\mathbb{R})\bigr),
$$
with
$$
\supp(w(\cdot,t))\cup\supp(F(\cdot,t))\subseteq[-1,1]
$$
for every $t\in[0,1)$, while the solution exhibits finite-time blow-up in the sense that
$$
\limsup_{t\nearrow1}\|w(\cdot,t)\|_{L^\infty}=\infty.
$$
\end{thm}

\subsection{Ideas of the proof}

We will build the solution announced in Theorem~\ref{construccionblowup} as an infinite superposition of smooth, odd, compactly supported functions, which we call \emph{layers}, switched on one after another at times accumulating at $t=1$. The vortex layer cascade method is inherently iterative so, for the sake of clarity, we first describe what happens with only two layers.

Let us start with a solution $w^1$ of Equation~\eqref{eq:DGF} with forcing $F^1$. We call $w^1$ the \emph{background}. This background must be odd in the space variable. We can perturb this solution by adding a \emph{perturbation} $w^2$. In that case, the sum $w^1+w^2$ solves Equation~\eqref{eq:DGF} with forcing $F^1+F^2$, where, writing $u^i=u(w^i)$,
$$
F^2=\underbrace{\partial_tw^2+u^1\partial_xw^2-w^2Hw^1}_{\textit{Evolution produced by the background}}+\underbrace{u^2\partial_xw^1-w^1Hw^2}_{\textit{Cross terms}}+\underbrace{u^2\partial_xw^2-w^2Hw^2}_{\textit{Quadratic terms}}.
$$

The whole point of the method is to choose $w^2$ with a specific functional form so that two things happen simultaneously. On the one hand, the forcing must remain regular: a rough forcing could produce a blow-up by itself, so the interest of the construction lies in making the forcing as regular as possible. Specifically, in this paper we obtain a forcing which is bounded in $C^{0,1-\delta}(\R)$ uniformly in time up to the blow-up time, for $\delta>0$ as small as desired.

On the other hand, the solution must grow: since Equation~\eqref{eq:DGF} is a toy model for the three-dimensional Euler equations in which $w$ plays the role of the vorticity, we would like the superposition of infinitely many layers to blow up in finite time.

Let us first look at the term describing the evolution of the perturbation under the background. Ideally, we would choose $w^2$ so that this term vanishes. Although the corresponding equation can be solved by characteristics, its solution involves the flow map of $u^1$, which makes the analysis cumbersome. We therefore approximate $u^1$ by its first-order Taylor polynomial at the origin. Since $u^1(0,t)=0$ and $\partial_xu^1=Hw^1$, this amounts to replacing $u^1(x,t)$ by $Hw^1(0,t)x$ and $Hw^1(x,t)$ by $Hw^1(0,t)$. We are thus led to the linearized equation
\begin{equation}\label{eq:lin_int}
    \partial_tw^2+Hw^1(0,t)x\partial_xw^2=Hw^1(0,t)w^2,
\end{equation}
which can be easily solved (see Lemma~\ref{DuhamelForm}).

Indeed, let $t_0$ be the time at which the perturbation is introduced. The characteristic curves of the linearized equation are the dilations
$$
X(t;t_0,x)=x\exp\left(\int_{t_0}^tHw^1(0,s)\dd s\right),
$$
and along them the perturbation satisfies
$$
\frac{\dd}{\dd t}w^2\bigl(X(t;t_0,x),t\bigr)=Hw^1(0,t)w^2\bigl(X(t;t_0,x),t\bigr),
$$
so that
$$
w^2\bigl(X(t;t_0,x),t\bigr)=\exp\left(\int_{t_0}^tHw^1(0,s)\dd s\right)w^2(x,t_0).
$$

Therefore, as long as $Hw^1(0,t)>0$, the background pushes the perturbation, spreading its support and amplifying it by the same exponential factor, so that a perturbation which is very small when it is introduced ends up being much larger than the background. This is the mechanism of growth of the vorticity on which the method relies. 

It will be proven in Theorem~\ref{thmIterativo} that if $Hw^1(0,t)\sim N^\varepsilon$ while the perturbation evolves, then $H(w^1+w^2)(0,t)\sim M^\varepsilon$ at the end of the evolution, where $M$ is exponentially larger than $N$.

The price to pay is that the error arising from this approximation of the velocity remains in the forcing term. For this error to be small, the Taylor approximation must be good enough, which forces $w^2$ to be much more concentrated near the origin than $w^1$.

The quadratic terms pose a different problem. They do not involve the background at all, and for an arbitrary perturbation there is no hope that they are small. The situation changes if the perturbation is highly oscillatory. More precisely, one could take $w^2(x,t_0)=\varphi(x)$ with
$$
\varphi(x)=-A\chi(\lambda x)\sin(kx),
$$
where $\chi$ is a smooth even cutoff, $A>0$ is a small amplitude and $k$ is much larger than $\lambda$. On such functions the Hilbert transform acts as
$$
H\varphi\approx-\frac{\varphi'}{k},\qquad \int_0^xH\varphi(y)\dd y\approx-\frac{\varphi}{k},
$$
with lower order errors.

Since the linearized evolution only dilates and amplifies the perturbation, $w^2(\cdot,t)$ has the same form at every time up to the different values of $A$, $\lambda$ and $k$. Consequently, the leading parts of $u^2\partial_xw^2$ and $w^2Hw^2$ coincide and cancel each other in the quadratic terms. This cancellation, which is made precise in Remark~\ref{cancelacionfund}, is the reason behind our choice of the profile.

\begin{figure}[htbp]
\centering
\pgfmathsetmacro{\mechB}{0}
\pgfmathsetmacro{\mechD}{4}
\pgfmathsetmacro{\mechAA}{2.1}
\pgfmathsetmacro{\mechSS}{0.17}
\pgfmathsetmacro{\mechA}{\mechAA/\mechD}
\pgfmathsetmacro{\mechS}{\mechSS/\mechD}
\begin{tikzpicture}
\begin{groupplot}[
  group style={group size=2 by 1, horizontal sep=1cm},
  width=5.6cm, height=4.2cm, scale only axis,
  xmin=-1.1, xmax=1.1, ymin=-2.3, ymax=2.3,
  axis lines=none, clip=false,
  title style={yshift=6pt, text height=2ex, text depth=0.5ex},
  every axis plot/.append style={line join=round, thick},
]
\nextgroupplot[title={Shortly after $w^2$ is switched on}]
\draw[-stealth] (axis cs:-1.1,0) -- (axis cs:1.1,0);
\addplot[blue, dash pattern=on 6pt off 3pt, line cap=round, domain=-1:1, samples=500] {bgw(x-\mechB*x*(1-x*x)*(1-x*x))};
\addplot[red, semithick, domain=-\mechS:\mechS, samples=300] {ptw(x,\mechA,\mechS)};
\draw[red,-latex] (axis cs:0,0.7) -- (axis cs:0,1.6);
\draw[red,-latex] (axis cs:0,-0.7) -- (axis cs:0,-1.6);
\draw[red,-latex] (axis cs:0.08,0.28) -- (axis cs:\mechSS,0.28);
\draw[red,-latex] (axis cs:-0.08,-0.28) -- (axis cs:-\mechSS,-0.28);
\node[blue] at (axis cs:-0.42,1.0) {$w^1$};
\node[red]  at (axis cs:0.22,0.75) {$w^2$};

\nextgroupplot[title={Almost at the end of the evolution}]
\draw[-stealth] (axis cs:-1.1,0) -- (axis cs:1.1,0);
\addplot[blue, dash pattern=on 6pt off 3pt, line cap=round, domain=-1:1, samples=500] {bgw(x-\mechB*x*(1-x*x)*(1-x*x))};
\addplot[red, domain=-\mechSS:\mechSS, samples=600] {ptw(x,\mechAA,\mechSS)};
\node[blue] at (axis cs:-0.42,1.0) {$w^1$};
\node[red]  at (axis cs:0.3,1.9) {$w^2$};
\end{groupplot}
\end{tikzpicture}
\caption{Schematic picture of the mechanism of growth with two layers. The perturbation $w^2$ is introduced near the origin with a very small amplitude and a very small support, compared with those of the background $w^1$. The background, which almost does not change during the perturbation evolution, dilates the perturbation, so that its amplitude and the width of its support grow by the same factor. Almost at the end of the evolution, $w^2$ is much larger than $w^1$, but still much more concentrated near the origin.}
\label{fig:mecanismo}
\end{figure}
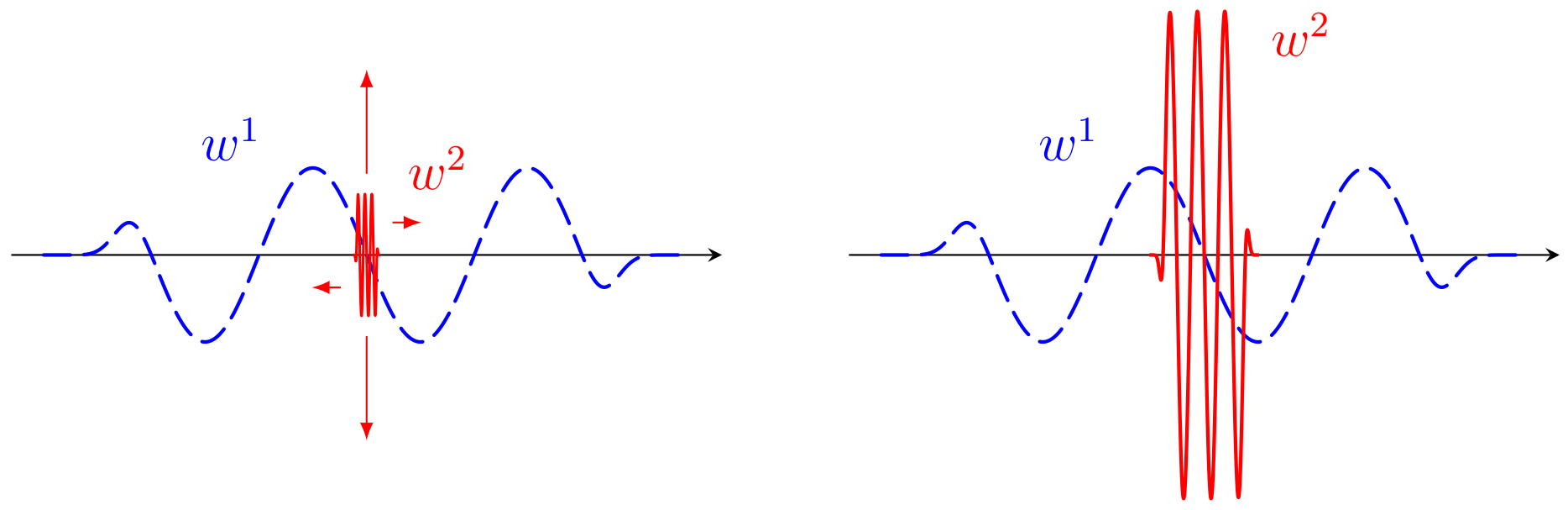

The strategy outlined above is adapted from the approach developed in~\cite{CordobaMartinezZoroaBlowupEuler}. When applied to Equation~\eqref{eq:DGF}, it yields a forcing term in $C^{0,1/2-\delta}$, as shown in~\cite{TFM}.

However, the method can be pushed further. As explained before, the perturbation has been chosen as a solution to Equation~\eqref{eq:lin_int}. Consequently, the remaining terms to be controlled are the errors arising from this approximation, together with the quadratic and cross terms.

The quadratic terms are unavoidable, since they stem from the nonlinearity of the equation, but the other two involve both the background and the perturbation in a very particular way. This raises the following possibility: if, instead of choosing the perturbation as a solution to Equation~\eqref{eq:lin_int}, we allowed it to solve the same equation with a suitably chosen forcing term designed to cancel or reduce some of these errors, could we obtain a forcing term $F^2$ with improved regularity?

The answer is yes, and this is part of what is done in~\cite{CordobaMartinezZoroaBlowupIPM} in order to obtain a smooth forcing. In our simpler setting we perform just one step of this procedure. We add to the perturbation two correctors which solve Equation~\eqref{eq:lin_int} with forcing terms chosen to cancel exactly the largest contributions from the linearization error and the cross terms. This raises the regularity of the forcing from $C^{0,1/2-\delta}$ to $C^{0,1-\delta}$.

In fact, the same idea can be iterated to make the linearization error and the cross terms as small as desired, at the expense of more complicated layers. This would not improve Theorem~\ref{construccionblowup}, which already holds for every $\delta\in(0,1)$. What prevents us from reaching Lipschitz regularity are the self-interaction of the leading profile, which is not affected by the correctors, and an activation term arising from switching on the layers through a cutoff function. We have chosen not to pursue this, since our aim is to present the vortex layer cascade method in its simplest form rather than to optimize the regularity of the forcing.

Finally, the method is iterated as follows: starting from a fixed background, we construct a perturbation following the procedure described above, and we then repeat the construction taking as the new background the sum of the previous background and its perturbation. Each new layer is introduced closer to $t=1$ and grows over a shorter time interval than the previous one, so that infinitely many layers are added before $t=1$ while the Hilbert transform of the superposition at the origin becomes unbounded. This leads to the blow-up of the solution, and the details are given in the proof of Theorem~\ref{construccionblowup}.

\subsection{Preliminaries and notation}\label{sec:prelim}

We write $A\lesssim B$ if $A\leq CB$ for some constant $C>0$ independent of all relevant parameters, and $A\lesssim_\gamma B$ when $C$ is allowed to depend on $\gamma$. If $A\lesssim B$ and $B\lesssim A$, we write $A\sim B$.

Throughout, $I\subset\R$ denotes a closed, possibly unbounded, interval. For $k\in\N_0$, $C^k(I)$ is the space of $k$ times continuously differentiable functions $f\colon I\to\R$ whose derivatives up to order $k$ are bounded, endowed with the norm $\|f\|_{C^k(I)}=\sum_{j=0}^k\|f^{(j)}\|_{L^\infty(I)}$. 

Accordingly, $C^\infty(I)=\bigcap_{k\in\N_0}C^k(I)$ is the space of smooth functions on $I$ all of whose derivatives are bounded. 

When we only know that a function is infinitely differentiable, with no control on the size of its derivatives, we say that it is \emph{smooth}. 

The same definitions apply to complex-valued functions and to functions of $(x,t)\in\R\times J$, with $J\subset \R$ a closed interval, by requiring all partial derivatives up to the corresponding order to be continuous and bounded on $\R\times J$.

For $b>0$ and $k\in\N_0\cup\{\infty\}$, we write $f\in C^k(\R\times[0,b))$ if $f\in C^k(\R\times[0,c])$ for every $c\in(0,b)$, in agreement with the notion of classical solution on $[0,b)$ given below.

Given $k\in\N_0$ and $\alpha\in(0,1)$, the Hölder space $C^{k,\alpha}(I)$ consists of those $f\in C^k(I)$ for which $\|f\|_{C^{k,\alpha}(I)}=\|f\|_{C^k(I)}+[f^{(k)}]_{C^{0,\alpha}(I)}<\infty$, where
$$
[g]_{C^{0,\alpha}(I)}=\sup_{\substack{x,y\in I\\ x\neq y}}\frac{|g(x)-g(y)|}{|x-y|^\alpha}.
$$
Whenever the domain is clear from the context, we omit it and simply write $C^k$, $C^{k,\alpha}$, $L^\infty$, $\|\cdot\|_{C^{k,\alpha}}$, and so on.

We will repeatedly use the following elementary inequality.
\begin{lem}[Interpolation]\label{interpolacion}
    Let $\alpha\in(0,1)$ and $f\in C^1(I)$. Then $f\in C^{0,\alpha}(I)$ and $[f]_{C^{0,\alpha}(I)}\leq 2^{1-\alpha}\|f\|_{L^\infty(I)}^{1-\alpha}\|f'\|_{L^\infty(I)}^{\alpha}$. In particular,
    \begin{equation}\label{des_interpolacion}
        \|f\|_{C^{0,\alpha}(I)}\leq \|f\|_{L^\infty(I)}+2^{1-\alpha}\|f\|_{L^\infty(I)}^{1-\alpha}\|f'\|_{L^\infty(I)}^{\alpha}.
    \end{equation}
\end{lem}
\begin{proof}
    Just interpolate between the bounds $|f(x)-f(y)|\leq 2\|f\|_{L^\infty(I)}$ and $|f(x)-f(y)|\leq\|f'\|_{L^\infty(I)}|x-y|$.
\end{proof}
Given a Banach space $X$ and $T>0$, we endow $C^k([0,T];X)$ with the norm $\|u\|_{C^k([0,T];X)}=\sum_{j=0}^k\sup_{t\in[0,T]}\|\partial_t^j u(t)\|_X$, and abbreviate it as $C^k_tX_x$ whenever the domains cannot be confused.

For $f\in\mathcal{S}(\R)$ (the Schwartz space on $\mathbb{R}$), its Hilbert transform is given by
$$
Hf(x)=\frac{1}{\pi}\operatorname{P.V.}\int_\R\frac{f(y)}{x-y}\dd y,
\qquad\text{so that}\qquad
\widehat{Hf}(\xi)=-i\operatorname{sgn}(\xi)\widehat{f}(\xi).
$$
We refer to~\cite{Grafakos2014} for the standard properties of the Hilbert transform used from now on, including its commutation with differentiation and positive dilations, and its boundedness on $L^p(\mathbb{R})$ for $1<p<\infty$. We will also use that the Hilbert transform reverses parity: if a given function $f$ is even, then $Hf$ is odd, and if $f$ is odd, then $Hf$ is even.

\subsubsection{Properties of the velocity}\label{sec:velocity}
In Equation~\eqref{eq:DGF} the velocity is recovered from $w$ through the one-dimensional analogue of the Biot--Savart law, $\partial_xu(w)=Hw$. As stated in the introduction, we will consider the following normalization
\begin{equation}\label{eq:BiotSavart1D}
u(w)(x)=\int_0^xHw(y)\dd y.
\end{equation}
In this subsection we collect the properties of this velocity operator that are used throughout the paper.
\begin{lem}[Hölder estimates]\label{lem:holderH}
Let $\beta\in(0,1)$, $k\in\N_0$ and $R>0$, and let $w\in C^{k,\beta}(\R)$ with $\supp(w)\subseteq[-R,R]$. Then $Hw\in C^{k,\beta}(\R)$ and the following estimates hold.
\begin{equation}\label{eq:privalov}
[Hw]_{C^{0,\beta}}\leq c(\beta)[w]_{C^{0,\beta}}
\end{equation}
for some $c(\beta)>0$,
\begin{equation}\label{eq:holderH}
\|Hw\|_{C^{k,\beta}}\lesssim_{c(\beta), \beta, k}(1+R)\|w\|_{C^{k,\beta}}.
\end{equation}
\end{lem}
\begin{proof}
    Estimate~\eqref{eq:privalov} is classical (see, for instance,~\cite{GarciaCuervaGatto2005} and~\cite{DyachenkoNursultanovTikhonov2013}). Inequality~\eqref{eq:holderH} follows from Estimate~\eqref{eq:privalov} by elementary arguments.
\end{proof}
Observe that, if $w\in C_c^{k,\beta}(\R)$ for some $k\in\N$, then $Hw\in C^k(\R)$ and therefore $u(w)$ is $k+1$ times continuously differentiable. Moreover, $\partial_x^{j+1}u(w)=H(\partial_x^jw)$ for $0\leq j\leq k$.

As a consequence of Lemma~\ref{lem:holderH} we also have the following result, whose proof follows by elementary considerations.
\begin{lem}[Boundedness of the velocity operator on Hölder spaces]\label{velocidadHolder}
		Let $0<\beta<1$ and $R>0$, and let $w\in C^{0, \beta}(\R)$ be an odd function with $\supp(w)\subseteq [-R, R].$ Then 
        $$\|u(w)\|_{L^\infty}\lesssim_{c(\beta), \beta} (1+R+R^2)\|w\|_{C^{0,\beta}}.$$
		In particular, if $k\in\N$ and, in addition, $w\in C^{k,\beta}(\R)$, then
		$$
		\|u(w)\|_{C^{k+1,\beta}}\lesssim_{c(\beta), \beta, k} (1+R+R^2)\|w\|_{C^{k,\beta}}.
		$$
	\end{lem}

The following logarithmic estimate will only be used at the end of the proof of Theorem~\ref{construccionblowup}, to deduce the blow-up of $\|w(\cdot,t)\|_{L^\infty}$ from the growth of $Hw(0,t)$.
\begin{lem}\label{estimalogH}
	Let $f\in C^1_c(\R)$ satisfy $f\not\equiv 0$ and $\supp(f)\subseteq[-R, R]$ for some $R\geq 1$. Then
	$$
	\|Hf\|_{L^\infty}\leq \frac{2}{\pi}\|f\|_{L^\infty}\left[2+\log\left(\frac{3R\|f\|_{C^1}}{\|f\|_{L^\infty}}\right)\right].
	$$
\end{lem}
We omit the proof since, again, it only relies on elementary arguments: one splits the principal value defining $Hf(x)$ into a near-field part where the cancellation of the kernel allows one to replace $f(x-y)$ by $f(x-y)-f(x)$ and to apply the mean value theorem, and a far-field part which is controlled using the support of $f$.

\subsubsection{Classical solutions}
Although the notion of a classical solution is fairly intuitive, the precise definition adopted throughout this paper is the following.

\begin{dfn}
Let $b>0$. We say that $w$ is a \emph{classical
solution for Equation~\eqref{eq:DGF} with forcing $F$} on $[0,b]$ if $w\in C^1(\mathbb R\times[0,b])$, $F\in C(\mathbb R\times[0,b])$ and the following conditions hold:
\begin{enumerate}
    \item The equation $\partial_t w+u(w)\partial_x w=wHw+F$ is satisfied pointwise in $\mathbb R\times(0,b)$, where
    $$
    u(w)(x,t)=\int_0^x Hw(y,t) \dd y, \quad \forall (x,t)\in \R\times[0,b]
    $$
    \item The initial condition $w(x,0)=w_0(x)$ is satisfied for every $x\in\R$.
    \item There exists a compact set $K\subset\mathbb R$, independent of time, such that
    $$
    \supp(w(\cdot,t))\cup \supp(F(\cdot,t))\subseteq K
    \qquad\text{for every }t\in[0,b].
    $$
\end{enumerate}
We say that $w$ is a classical solution with forcing $F$ on $[0,b)$ if its restriction to $[0,c]$ is a classical solution on $[0,c]$ for every $c\in(0,b)$, and there exists a single compact set $K\subset\mathbb R$ such that
    $$
    \supp(w(\cdot,t))\cup \supp(F(\cdot,t))\subseteq K
    \qquad\text{for every }t\in[0,b).
    $$
\end{dfn}
\begin{rem}
    Whenever the compact set is left unspecified, its existence is implicit, and its particular choice is irrelevant to the discussion.
\end{rem}

\section{Approximation estimates for the Hilbert transform and its primitive}\label{seccion2}

Throughout this paper $\chi$ will denote a nonnegative, even bump function in $C_c^\infty(\mathbb{R})$, supported in $[-1,1]$ and identically equal to $1$ on $[-1/2,1/2]$. Moreover, $0\leq\chi(x)\leq 1$ for every $x\in\R$.
\subsection{General estimates}

\begin{lem}\label{auxiliarwave}
Let $n,m\in\mathbb{N}_0$ and $0<\lambda<k<\infty$, and define
$$
g^{(n,m)}_{\lambda,k}(x)=x^n\chi^{(m)}(\lambda x)e^{ikx}.
$$
Then there exists $R^{(n,m)}_{\lambda,k}\in C^\infty(\mathbb{R})$ such that
\begin{equation}\label{eq:auxiliarwave}
    Hg^{(n,m)}_{\lambda,k}=-ig^{(n,m)}_{\lambda,k}+R^{(n,m)}_{\lambda,k},
\end{equation}
with
$$
\|\partial_x^j R^{(n,m)}_{\lambda,k}\|_{L^\infty}\leq C^{(n,m)}_{\chi,K,j}\frac{\lambda^{j-n+K}}{k^K}
$$
for every $j\in\mathbb{N}_0$ and $K\in\mathbb{N}$, where $C^{(n,m)}_{\chi,K,j}>0$ is a constant independent of $\lambda$ and $k$.
\end{lem}

\begin{proof}
Since $g^{(n,m)}_{\lambda,k}\in\mathcal{S}(\mathbb{R})$ and the
Hilbert transform is bounded on $H^s(\mathbb{R})$ for every $s\geq 0$,
we have
$$
R^{(n,m)}_{\lambda,k}=Hg^{(n,m)}_{\lambda,k}+ig^{(n,m)}_{\lambda,k}\in H^s(\mathbb{R})
$$
for every $s\geq0$. In particular, $R^{(n,m)}_{\lambda,k}\in C^\infty(\mathbb{R})$ and $\partial_x^jR^{(n,m)}_{\lambda,k}\in L^2(\mathbb{R})$ for every $j\in\mathbb{N}_0$.

Fix $j\in\mathbb{N}_0$. By the Fourier multiplier formula for the Hilbert transform, we have
$$
\left[\partial_x^j R^{(n,m)}_{\lambda,k}\right]^\wedge(\xi)=(i\xi)^j\left(1-\operatorname{sgn}(\xi)\right)i\widehat{g^{(n,m)}_{\lambda,k}}(\xi)=2i\ind_{(-\infty,0)}(\xi)(i\xi)^j\widehat{g^{(n,m)}_{\lambda,k}}(\xi)
$$
for almost every $\xi\in\mathbb{R}$. Since
$g^{(n,m)}_{\lambda,k}\in\mathcal{S}(\mathbb{R})$, it follows that $\left[\partial_x^j R^{(n,m)}_{\lambda,k}\right]^\wedge
\in L^1(\mathbb{R}).$
For $a>0$, define
$$
F_a(x)=\frac{1}{2\pi}\int_{-a}^{a}
\left[\partial_x^j R^{(n,m)}_{\lambda,k}\right]^\wedge(\xi)
e^{ix\xi}\dd\xi.
$$
By the Fourier inversion theorem in $L^2(\mathbb{R})$ and
Plancherel's theorem, $F_a$ converges to
$\partial_x^j R^{(n,m)}_{\lambda,k}$ in $L^2(\mathbb{R})$ as
$a\to\infty$. Hence, there exists a sequence $a_\ell\to\infty$
such that $F_{a_\ell}(x)\to
\partial_x^j R^{(n,m)}_{\lambda,k}(x)$
for almost every $x\in\mathbb{R}$. Consequently,
$$
\left|\partial_x^j R^{(n,m)}_{\lambda,k}(x)\right|\leq \frac{1}{2\pi}
\int_{\mathbb{R}}
\left|\left[\partial_x^j R^{(n,m)}_{\lambda,k}\right]^\wedge(\xi)\right|
\dd\xi=\frac{1}{\pi}\int_{-\infty}^0
|\xi|^j\left|\widehat{g^{(n,m)}_{\lambda,k}}(\xi)\right|\dd\xi
$$
for almost every $x\in\mathbb{R}$. Since
$\partial_x^j R^{(n,m)}_{\lambda,k}$ is continuous, the same
inequality holds for every $x\in\mathbb{R}$.

Now, let $G^{(n,m)}(x)=x^n\chi^{(m)}(x).$ Then $G^{(n,m)}\in\mathcal{S}(\mathbb{R})$ and
$g^{(n,m)}_{\lambda,k}(x)
=\lambda^{-n}G^{(n,m)}(\lambda x)e^{ikx}.$
Therefore,
$$
\widehat{g^{(n,m)}_{\lambda,k}}(\xi)
=\lambda^{-n-1}
\widehat{G^{(n,m)}}\left(\frac{\xi-k}{\lambda}\right).
$$
It follows that
$$
\left|\partial_x^j R^{(n,m)}_{\lambda,k}(x)\right|\leq \frac{1}{\pi}
\int_{-\infty}^0
\frac{|\xi|^j}{\lambda^{n+1}}
\left|\widehat{G^{(n,m)}}\left(\frac{\xi-k}{\lambda}\right)\right|
\dd\xi=\frac{1}{\pi}\lambda^{-n}
\int_{-\infty}^{-k/\lambda}
|k+\lambda u|^j\left|\widehat{G^{(n,m)}}(u)\right|\dd u.
$$
Since $0<\lambda<k$, for $u\leq-k/\lambda$ we have $|k+\lambda u|=-k-\lambda u\leq-\lambda u=\lambda|u|.$
Consequently,
$$
\left|\partial_x^j R^{(n,m)}_{\lambda,k}(x)\right|\leq \frac{1}{\pi}\lambda^{j-n}
\int_{-\infty}^{-k/\lambda}
|u|^j\left|\widehat{G^{(n,m)}}(u)\right|\dd u.
$$
For any $K\in\mathbb{N}$, the fact that
$G^{(n,m)}\in\mathcal{S}(\mathbb{R})$ implies that
$$
c^{(n,m)}_{\chi,K,j}=\frac{1}{\pi}\sup_{u\in\mathbb{R}}|u|^{j+K+1}
\left|\widehat{G^{(n,m)}}(u)\right|<\infty.
$$
Thus,
$$
\left|\partial_x^j R^{(n,m)}_{\lambda,k}(x)\right|\leq c^{(n,m)}_{\chi,K,j}\lambda^{j-n}
\int_{k/\lambda}^{\infty}\frac{du}{u^{K+1}}=\frac{c^{(n,m)}_{\chi,K,j}}{K}\frac{\lambda^{j-n+K}}{k^K}=C^{(n,m)}_{\chi,K,j}\frac{\lambda^{j-n+K}}{k^K}.
$$
Taking the supremum over $x\in\mathbb{R}$ yields the desired
estimate.
\end{proof}

\begin{lem}\label{auxiliarwaveP}
Let $n,m\in\mathbb{N}_0$ and $0<\lambda<k<\infty$, and let
$g^{(n,m)}_{\lambda,k}(x)=x^n\chi^{(m)}(\lambda x)e^{ikx}$ be as in Lemma~\ref{auxiliarwave}.
Then, for every $K\in\mathbb{N}$, there exists a smooth function $s^{(n,m)}_{K,\lambda,k}$ such that, for all $x\in\mathbb{R}$,
\begin{equation}\label{eq:auxiliarwaveP}
    \int_0^xHg^{(n,m)}_{\lambda,k}(y)\dd y
=-i\sum_{r=0}^{K-1}\frac{(-1)^r}{(ik)^{r+1}}\left[
\partial_x^r\left(x^n\chi^{(m)}(\lambda x)\right)e^{ikx}
-
\partial_x^r\left(x^n\chi^{(m)}(\lambda x)\right)\Big|_{x=0}\right]+s^{(n,m)}_{K,\lambda,k}(x).
\end{equation}
Moreover, there exists a constant $C^{(n,m)}_{\chi,K}>0$, independent of $\lambda$, $k$ and $L$, such that for every $L>0$,
$$
\|s^{(n,m)}_{K,\lambda,k}\|_{L^\infty([-L,L])}
\leq
C^{(n,m)}_{\chi,K}L\frac{\lambda^{K-n}}{k^K},
\qquad
\|\partial_xs^{(n,m)}_{K,\lambda,k}\|_{L^\infty(\mathbb{R})}
\leq
C^{(n,m)}_{\chi,K}\frac{\lambda^{K-n}}{k^K}.
$$
\end{lem}

\begin{proof}
Let $a^{(n,m)}_\lambda(x)=x^n\chi^{(m)}(\lambda x)$, so that
$g^{(n,m)}_{\lambda,k}(x)=a^{(n,m)}_\lambda(x)e^{ikx}.$
By Lemma~\ref{auxiliarwave}, $Hg^{(n,m)}_{\lambda,k}=-ig^{(n,m)}_{\lambda,k}+R^{(n,m)}_{\lambda,k}$, and hence
$$
\int_0^xHg^{(n,m)}_{\lambda,k}(y)\dd y=-i\int_0^xa^{(n,m)}_\lambda(y)e^{iky}\dd y+\int_0^xR^{(n,m)}_{\lambda,k}(y)\dd y.
$$
Integrating by parts once gives
$$
\int_0^xa^{(n,m)}_\lambda(y)e^{iky}\dd y=\frac{a^{(n,m)}_\lambda(x)e^{ikx}-a^{(n,m)}_\lambda(0)}{ik}-\frac{1}{ik}\int_0^x\partial_ya^{(n,m)}_\lambda(y)e^{iky}\dd y.
$$
Applying this identity repeatedly, we obtain
\begin{align*}
\int_0^xa^{(n,m)}_\lambda(y)e^{iky}\dd y&=\sum_{r=0}^{K-1}\frac{(-1)^r}{(ik)^{r+1}}\left[\partial_x^ra^{(n,m)}_\lambda(x)e^{ikx}-\partial_x^ra^{(n,m)}_\lambda(0)\right]\\
&+\frac{(-1)^K}{(ik)^K}\int_0^x\partial_y^Ka^{(n,m)}_\lambda(y)e^{iky}\dd y.
\end{align*}
Multiplying by $-i$ and substituting into the previous expression, the claimed identity holds with
$$
s^{(n,m)}_{K,\lambda,k}(x)=-i\frac{(-1)^K}{(ik)^K}\int_0^x\partial_y^Ka^{(n,m)}_\lambda(y)e^{iky}\dd y+\int_0^xR^{(n,m)}_{\lambda,k}(y)\dd y.
$$
Clearly, $s^{(n,m)}_{K,\lambda,k}$ is smooth. It remains to prove the estimates. Let $G^{(n,m)}(x)=x^n\chi^{(m)}(x)$ as in Lemma~\ref{auxiliarwave}. Then $G^{(n,m)}\in C_c^\infty(\mathbb{R})$ and
$a^{(n,m)}_\lambda(x)=\lambda^{-n}G^{(n,m)}(\lambda x)$. Therefore
$\partial_x^Ka^{(n,m)}_\lambda(x)=\lambda^{K-n}\left(G^{(n,m)}\right)^{(K)}(\lambda x)$, so that
$$
\|\partial_x^Ka^{(n,m)}_\lambda\|_{L^\infty(\mathbb{R})}\leq\lambda^{K-n}\|(G^{(n,m)})^{(K)}\|_{L^\infty(\mathbb{R})}.
$$
Recall that, by Lemma~\ref{auxiliarwave} with $j=0$, we also have that
$$
\|R^{(n,m)}_{\lambda,k}\|_{L^\infty(\mathbb{R})}\leq C^{(n,m)}_{\chi,K,0}\frac{\lambda^{K-n}}{k^K}.
$$
Set $C^{(n,m)}_{\chi,K}=\|\partial^K_xG^{(n,m)}\|_{L^\infty(\mathbb{R})}+C^{(n,m)}_{\chi,K,0}$ and let $L>0$ and $x\in[-L,L]$. We then get
$$
\left|s^{(n,m)}_{K,\lambda,k}(x)\right|
\leq|x|\left(\frac{1}{k^K}\|\partial_x^Ka^{(n,m)}_\lambda\|_{L^\infty(\mathbb{R})}+\|R^{(n,m)}_{\lambda,k}\|_{L^\infty(\mathbb{R})}\right)
\leq C^{(n,m)}_{\chi,K}L\frac{\lambda^{K-n}}{k^K}.
$$
Taking the supremum over $x\in[-L,L]$ gives the first estimate.

Finally, differentiating the definition of $s^{(n,m)}_{K,\lambda,k}$ we obtain
$$
\partial_xs^{(n,m)}_{K,\lambda,k}(x)
=-i\frac{(-1)^K}{(ik)^K}\partial_x^Ka^{(n,m)}_\lambda(x)e^{ikx}+R^{(n,m)}_{\lambda,k}(x),
$$
and therefore
$$
\|\partial_xs^{(n,m)}_{K,\lambda,k}\|_{L^\infty(\mathbb{R})}\leq\frac{1}{k^K}\|\partial_x^Ka^{(n,m)}_\lambda\|_{L^\infty(\mathbb{R})}+\|R^{(n,m)}_{\lambda,k}\|_{L^\infty(\mathbb{R})}\leq C^{(n,m)}_{\chi,K}\frac{\lambda^{K-n}}{k^K}.
$$
This concludes the proof.
\end{proof}

\begin{lem}\label{lem:estima_dj_gen}
Let $n,m\in\mathbb{N}_0$ and $0<\lambda<k<\infty$, and let
$g^{(n,m)}_{\lambda,k}(x)=x^n\chi^{(m)}(\lambda x)e^{ikx}$ be as in Lemma~\ref{auxiliarwave}.
Then for every $j\in\N_0$ we have that
\begin{equation}\label{eq:estima_dj_gen}
    \|\partial^j_x g^{(n,m)}_{\lambda,k}\|_{L^\infty}\leq \sum_{\ell=0}^{\min\{n,j\}}\binom{j}{\ell}\frac{n! 2^{j-\ell}}{(n-\ell)!}\|\chi\|_{C^{m+j-\ell}}\lambda^{\ell-n}k^{j-\ell}.
\end{equation}
In particular, 
\begin{equation}\label{eq:estima_dj_util}
    \|\partial^j_x g^{(n,m)}_{\lambda, k}\|_{L^\infty}\lesssim_{\chi, n, m, j}k^j\lambda^{-n}.
\end{equation}
\end{lem}
\begin{proof}
Let $x\in\R$ and suppose first that $n=0$. By Leibniz's rule,
$$
\partial^j_xg^{(0,m)}_{\lambda,k}(x)=\sum_{\ell=0}^j\binom{j}{\ell}\lambda^\ell\chi^{(m+\ell)}(\lambda x)(ik)^{j-\ell}e^{ikx}.
$$
Since $\lambda<k$, we have $\lambda^\ell k^{j-\ell}\leq k^j$ for every $0\leq\ell\leq j$, and therefore
$$
\|\partial^j_xg^{(0,m)}_{\lambda,k}\|_{L^\infty}\leq k^j\sum_{\ell=0}^j\binom{j}{\ell}\|\chi^{(m+\ell)}\|_{L^\infty}\leq 2^jk^j\|\chi\|_{C^{m+j}}.
$$
Suppose now that $n\in\N$. Applying Leibniz's rule again,
\begin{align*}
        \partial^j_xg^{(n,m)}_{\lambda, k}(x)&=\partial^j_x\left(x^ng^{(0,m)}_{\lambda, k}(x)\right)=\sum_{\ell=0}^j\binom{j}{\ell}\partial^\ell_x\left(x^n\right)\partial^{j-\ell}_xg^{(0,m)}_{\lambda, k}(x)\\&=\sum_{\ell=0}^{\min\{n,j\}}\binom{j}{\ell}\frac{n!}{(n-\ell)!}x^{n-\ell}\partial^{j-\ell}_xg^{(0,m)}_{\lambda, k}(x).
    \end{align*}
Since $\supp(\chi)\subseteq[-1,1]$, every derivative of $g^{(0,m)}_{\lambda,k}$ vanishes unless $|x|\leq\lambda^{-1}$, in which case $|x|^{n-\ell}\leq\lambda^{\ell-n}$. Estimate~\eqref{eq:estima_dj_gen} then follows from the case $n=0$ applied with $j-\ell$ in place of $j$. Regarding Estimate~\eqref{eq:estima_dj_util}, note that it follows from Estimate~\eqref{eq:estima_dj_gen} because $\lambda^\ell k^{-\ell}\leq 1$.
\end{proof}

\subsection{Building blocks of the layers}

Although the approximation results established in the previous subsection hold to arbitrary order and for a broad class of profiles, we now specialize them to the particular functions that will be used in the proof of Theorem~\ref{construccionblowup} to construct the singular solution.

\begin{dfn}\label{defperfiles}
		Let $M> 1$ and $\chi$ be as at the beginning of Section~\ref{seccion2}. For each $x\in\R$ define:
		\begin{itemize}
			\item $\varphi_0(x)=-\operatorname{Im}\left(g^{(0,0)}_{M^{1/2},M}(x)\right)=-\chi(M^{1/2}x)\sin(Mx)$,
			\item $\begin{aligned}[t]
				\varphi_1(x)&=x^3\varphi_0'(x)=-M^{1/2}\operatorname{Im}\left(g^{(3,1)}_{M^{1/2},M}(x)\right)-M\operatorname{Re}\left(g^{(3,0)}_{M^{1/2},M}(x)\right)\\
				&=-M^{1/2}x^3\chi'(M^{1/2}x)\sin(Mx)-Mx^3\chi(M^{1/2}x)\cos(Mx),
			\end{aligned}$
			\item $\varphi_2(x)=\operatorname{Re}\left(g^{(1,0)}_{M^{1/2},M}(x)\right)=x\chi(M^{1/2}x)\cos(Mx)$.
		\end{itemize}
\end{dfn}

\begin{rem}\label{rem:perfiles}
    Observe that $\varphi_i\in C^\infty_c(\R)$ is odd and $\supp(\varphi_i)\subseteq\{|x|\leq M^{-1/2}\}$ for each $i\in\{0,1,2\}$.
\end{rem}    

The next corollary collects the approximations of $H\varphi_i$ and of its primitive that will be used throughout. The key point is that both are, to leading order, multiples of $\varphi_i'$ and $\varphi_i$, respectively.

\begin{cor}\label{lem:aglutinado}
With the notation of Definition~\ref{defperfiles}, for each $i\in\{0,1,2\}$ there exist smooth functions $r_i$ and $s_i$, and $e_i, E_i\in C_c^\infty(\R)$, the latter two supported in ${\{|x|\leq M^{-1/2}\}}$, such that for every $x\in\R$ the following identities hold:
\begin{align}
    H\varphi_i(x)&=-\frac{\varphi_i'(x)}{M}+e_i(x)+r_i(x), \label{reescrituraperfiles} \\
    \int_0^xH\varphi_i(y)\dd y&=-\frac{\varphi_i(x)}{M}+E_i(x)+s_i(x). \label{reescrituraperfilesP}
\end{align}
In addition, for every $x\in\R$,
\begin{align}
    H\varphi_0(x)&=\chi(M^{1/2}x)\cos(Mx)+r_0(x), \label{Hphi0}\\
    E_0(x)&=\frac{\chi'(M^{1/2}x)\cos(Mx)}{M^{3/2}}-\frac{\chi''(M^{1/2}x)\sin(Mx)}{M^{2}}, \label{E0explicita}
\end{align}
and for $i\in\{1,2\}$,
\begin{equation}\label{Hphiorigen}
    H\varphi_i(0)=r_i(0).
\end{equation}
Moreover, for every $j\in\N_0$, $K\in\N$, $L>0$, $l\in\{0,1,2\}$ and $i\in\{1,2\}$,
\begin{align}
&\|\partial_x^jr_0\|_{L^\infty}\lesssim_{\chi,K,j}M^{\frac{j-K}{2}},
&&\|\partial_x^jr_i\|_{L^\infty}\lesssim_{\chi,K,j}M^{\frac{j-1-K}{2}}, \label{cotasr}\\
&\|\partial_x^je_0\|_{L^\infty}\lesssim_{\chi,j}M^{j-\frac12},
&&\|\partial_x^je_i\|_{L^\infty}\lesssim_{\chi,j}M^{j-1}, \label{cotase}\\
&\|\partial_x^jE_0\|_{L^\infty}\lesssim_{\chi,j}M^{j-\frac32},
&&\|\partial_x^jE_i\|_{L^\infty}\lesssim_{\chi,j}M^{j-2}, \label{cotasE}\\
&\|s_l\|_{L^\infty([-L,L])}\lesssim_{\chi}M^{-3/2}L,
&&\|\partial_xs_l\|_{L^\infty}\lesssim_{\chi}M^{-3/2}. \label{cotass}
\end{align}
\end{cor}

\begin{proof}
Throughout, $\lambda=M^{1/2}$ and $k=M$, and we write $g^{(n,m)}$, $R^{(n,m)}$ and $s^{(n,m)}_{K}$ for the functions of Lemmas~\ref{auxiliarwave} and~\ref{auxiliarwaveP} with these parameters, so that
\begin{equation}\label{eq:cotaR}
\|\partial_x^jR^{(n,m)}\|_{L^\infty}\lesssim_{n,m,\chi,K,j}M^{\frac{j-n-K}{2}},
\end{equation}
\begin{equation}\label{eq:cotaS}
    \|s^{(n,m)}_{K}\|_{L^\infty([-L,L])}\lesssim_{n,m,\chi,K}M^{-\frac{n+K}{2}}L, \qquad \|\partial_xs^{(n,m)}_{K}\|_{L^\infty}\lesssim_{n,m,\chi,K}M^{-\frac{n+K}{2}}.
\end{equation}
Recall also that the Hilbert transform maps real functions to real functions, so it commutes with taking real and imaginary parts. Taking real and imaginary parts in Equation~\eqref{eq:auxiliarwave} we obtain Identity~\eqref{Hphi0} together with
\begin{align*}
H\varphi_1(x)&=M^{1/2}x^3\chi'(M^{1/2}x)\cos(Mx)-Mx^3\chi(M^{1/2}x)\sin(Mx)+r_1(x),\\
H\varphi_2(x)&=x\chi(M^{1/2}x)\sin(Mx)+r_2(x),
\end{align*}
where
$$
r_0=-\operatorname{Im}\left(R^{(0,0)}\right),\qquad
r_1=-M^{1/2}\operatorname{Im}\left(R^{(3,1)}\right)-M\operatorname{Re}\left(R^{(3,0)}\right),\qquad
r_2=\operatorname{Re}\left(R^{(1,0)}\right).
$$
Evaluating these last two identities at $x=0$ gives Identity~\eqref{Hphiorigen}, and Estimates~\eqref{cotasr} follow from Estimate~\eqref{eq:cotaR}. We now define $e_i=H\varphi_i-r_i+\varphi_i'/M$, so that Identity~\eqref{reescrituraperfiles} holds by definition. Computing $\varphi_i'$ from Definition~\ref{defperfiles}, we find
\begin{align*}
e_0(x)&=-\frac{\chi'(M^{1/2}x)\sin(Mx)}{M^{1/2}},\\
e_1(x)&=-M^{1/2}x^3\chi'(M^{1/2}x)\cos(Mx)-x^3\chi''(M^{1/2}x)\sin(Mx)\\
&\quad-3x^2\chi(M^{1/2}x)\cos(Mx)-\frac{3x^2\chi'(M^{1/2}x)\sin(Mx)}{M^{1/2}},\\
e_2(x)&=\frac{\chi(M^{1/2}x)\cos(Mx)}{M}+\frac{x\chi'(M^{1/2}x)\cos(Mx)}{M^{1/2}}.
\end{align*}
We next turn to the primitives. We apply Lemma~\ref{auxiliarwaveP} with $K=3$ to $g^{(0,0)}$ and with $K=2$ to $g^{(3,1)}$, $g^{(3,0)}$ and $g^{(1,0)}$. Writing $a(x)=x^n\chi^{(m)}(\lambda x)$ and using that $a^{(r)}(0)=0$ for $r<n$ and that $\chi^{(r)}(0)=0$ for $r\geq1$, Identity~\eqref{eq:auxiliarwaveP} gives
\begin{align*}
\int_0^xHg^{(0,0)}(y)\dd y&=-\frac{a(x)e^{ikx}-1}{k}-\frac{ia'(x)e^{ikx}}{k^2}+\frac{a''(x)e^{ikx}}{k^3}+s^{(0,0)}_3(x),\\
\int_0^xHg^{(n,m)}(y)\dd y&=-\frac{a(x)e^{ikx}}{k}-\frac{i\bigl(a'(x)e^{ikx}-a'(0)\bigr)}{k^2}+s^{(n,m)}_2(x)\qquad(n\geq1).
\end{align*}
Taking real and imaginary parts as before and subtracting $-\varphi_i/M$, we obtain~\eqref{reescrituraperfilesP} with $E_0$ as in~\eqref{E0explicita},
\begin{align*}
E_1(x)&=-\frac{3x^2\chi(M^{1/2}x)\sin(Mx)}{M}-\frac{x^3\chi'(M^{1/2}x)\sin(Mx)}{M^{1/2}}\\
&\quad+\frac{3x^2\chi'(M^{1/2}x)\cos(Mx)}{M^{3/2}}+\frac{x^3\chi''(M^{1/2}x)\cos(Mx)}{M},\\
E_2(x)&=\frac{\chi(M^{1/2}x)\sin(Mx)}{M^2}+\frac{x\chi'(M^{1/2}x)\sin(Mx)}{M^{3/2}},
\end{align*}
and
$$
s_0=-\operatorname{Im}\left(s^{(0,0)}_3\right),\qquad
s_1=-M^{1/2}\operatorname{Im}\left(s^{(3,1)}_2\right)-M\operatorname{Re}\left(s^{(3,0)}_2\right),\qquad
s_2=\operatorname{Re}\left(s^{(1,0)}_2\right).
$$
Estimates~\eqref{cotass} then follow from~\eqref{eq:cotaS}. It remains to prove the assertions about $e_i$ and $E_i$. Each of these functions is a linear combination of real and imaginary parts of functions $g^{(n,m)}$, so it is supported in $\{|x|\leq M^{-1/2}\}$. Moreover, Bounds~\eqref{cotase} and~\eqref{cotasE} follow from Estimate~\eqref{eq:estima_dj_util}, which in this context gives $\|\partial_x^jg^{(n,m)}\|_{L^\infty}\lesssim M^{j-n/2}$.
\end{proof}

\section{The iterative result}
In this section, we prove the key step in the proof of Theorem~\ref{construccionblowup}, which will allow us to construct the singular solution iteratively.

\begin{thm}\label{thmIterativo}
Let $0<\varepsilon<1/9$ and let $N\geq N(\varepsilon, \chi)$ for some $N(\varepsilon, \chi)> 1$ sufficiently large. Suppose that $w^1\in C^\infty(\mathbb{R}\times[0,1])$ solves Equation~\eqref{eq:DGF} with forcing term $F^1\in C^\infty(\mathbb{R}\times[0,1])$. Assume that $w^1(\cdot,t)$ is odd for every $t\in[0,1]$, and that ${F^1\in C\bigl([0,1];C^{0,1-9\varepsilon}(\mathbb{R})\bigr).}$
Moreover, suppose that ${\supp\bigl(w^1(\cdot,t)\bigr)\cup \supp\bigl(F^1(\cdot,t)\bigr)\subseteq[-1,1]}$ for every $t\in[0,1]$ and, for each $t\in[1-N^{-\varepsilon/2},1]$,
$$
\|w^1(\cdot,t)\|_{C^{4,1/2}}\leq N^6,
$$
and
\begin{equation}\label{hipotesisTHen0}
\frac{N^\varepsilon}{4}\leq Hw^1(0,t)\leq 2N^\varepsilon.
\end{equation}
Define $M$ by
$$
M^{1/\varepsilon}=\exp\left(
\int_{1-N^{-\varepsilon/2}}^1 Hw^1(0,s)\dd s
\right).
$$

Then there exist $w^2\in C^\infty(\mathbb{R}\times[0,1])$ and a forcing term $F^2\in C^\infty(\mathbb{R}\times[0,1])$ such that $w^2(\cdot,t)$ is odd for every $t\in[0,1]$, $F^2\in C\bigl([0,1];C^{0,1-9\varepsilon}(\mathbb{R})\bigr)$ and \(w^1+w^2\) solves Equation~\eqref{eq:DGF} with forcing term \(F^1+F^2\). Furthermore, the forcing term \(F^2\) satisfies
\begin{equation}\label{eq:estimafuerza}
    \sup_{t\in[0,1]}\|F^2(\cdot,t)\|_{C^{0,1-9\varepsilon}}
\leq C_{\chi,\varepsilon}M^{-2\varepsilon},
\end{equation}
where \(C_{\chi,\varepsilon}>0\) is a constant depending only on \(\chi\) and \(\varepsilon\).

In addition, we have that $\supp\bigl((w^1+w^2)(\cdot,t)\bigr)\cup \supp\bigl((F^1+F^2)(\cdot,t)\bigr)\subseteq[-1,1]$ for every $t\in[0,1]$, and for each \(t\in[1-M^{-\varepsilon/2},1]\), the following estimates hold:
\begin{equation}\label{eq:reg_pert}
    \|(w^1+w^2)(\cdot,t)\|_{C^{4,1/2}}
\leq M^6,
\end{equation}
and
\begin{equation}\label{HIPOTESISTHEN0}
\frac{M^\varepsilon}{4}\leq H(w^1+w^2)(0,t)\leq 2M^\varepsilon.
\end{equation}
Finally,  for every \(0\leq t\leq 1-N^{-\varepsilon/2}\) we have that
\begin{equation}\label{eq:antesdetiempo}
    w^2(\cdot,t)=F^2(\cdot,t)=0.
\end{equation}
\end{thm}

\begin{rem}[Comparison between $N$ and $M$]\label{RelNM}
Throughout this section, we will repeatedly use the following
elementary comparison: for any $A,B,c,d>0$, there exists
$N_0=N_0(A,B,c,d,\varepsilon)$ such that $AN^c\leq BM^d$ for all $N\geq N_0$.

Indeed, by \eqref{hipotesisTHen0} and the definition of $M$, we have
$$
\log (M)
=\varepsilon\int_{1-N^{-\varepsilon/2}}^1 Hw^1(0,s)\dd s
\geq \frac{\varepsilon N^{\varepsilon/2}}{4}.
$$
Consequently, $M^d\geq e^{\frac{d\varepsilon}{4}N^{\varepsilon/2}}.$ Since $N^ce^{-\frac{d\varepsilon}{4}N^{\varepsilon/2}}
\to 0$ as $N\to\infty$, there exists $N_0=N_0(A,B,c,d,\varepsilon)$ such that
$$
N^c\exp\left(-\frac{d\varepsilon}{4}N^{\varepsilon/2}\right)
\leq \frac{B}{A}
\qquad\text{for all }N\geq N_0.
$$
It follows that
$$
AN^c
\leq B\exp\left(\frac{d\varepsilon}{4}N^{\varepsilon/2}\right)
\leq BM^d
\qquad\text{for all }N\geq N_0,
$$
as claimed.

Since only finitely many such comparisons will be needed
throughout this section, we may henceforth assume, without
further mention, that $N(\varepsilon, \chi)$ is sufficiently large for all of
them to hold simultaneously.
\end{rem}

\begin{proof}[Proof of Theorem~\ref{thmIterativo}]
The proof is divided into several steps, each of which is developed in a separate subsection below. Throughout, we refer to $w^1$ as the \emph{background layer} and to $w^2$ as the \emph{perturbation}. 

In Subsection~\ref{step1} we construct the perturbation and the associated forcing term, whose precise functional form will be essential to exploit the key cancellation in the quadratic term. In Subsections~\ref{step2} and~\ref{step3} we show that the perturbed layer $w^1+w^2$ satisfies all the desired properties. Finally, in Subsection~\ref{step4} we establish the required Hölder estimate on the forcing term.

For the sake of brevity, all statements in this section are understood to hold under the hypotheses of Theorem~\ref{thmIterativo}, and we shall freely use the notation introduced therein without further comment. 

\subsection{Construction of the perturbation and the forcing term}\label{step1}
For each $t\in[0,1]$ define
$$
a(t)=Hw^1(0,t), \qquad K(t)=\exp\left(\int_t^1 a(s)\dd s\right),
\qquad
b(t)=\frac{\log (K(t))}{\log (M)},
$$ 
and consider the transport operator $L=\partial_t+a(t)x\partial_x-a(t),$ that is,
$$
Lf(x,t)=\partial_tf(x,t)+a(t)x\partial_xf(x,t)-a(t)f(x,t).
$$
\begin{lem}[Duhamel formula for $L$]\label{DuhamelForm}
	Let $h\in C^\infty(\mathbb{R}\times[0,1])$. Then the function
	$$
	f(x,t)=\int_{1-N^{-\varepsilon/2}}^tM^{b(s)-b(t)}h\left(M^{b(t)-b(s)}x,s\right)\dd s
	$$
	is smooth and satisfies $Lf=h.$
	\end{lem}
	\begin{proof}
	The regularity of $f$ follows from the regularity of $h$ and from the fact that the functions $a$, $K$ and $b$ defined above are smooth. To verify the identity, we look for a solution of the form
	$$
	f(x,t)=M^{\varepsilon-b(t)}
	\psi\left(M^{b(t)-\varepsilon}x,t\right).
	$$
	Since $M^{b(t)}=K(t)$ and $K'(t)=-a(t)K(t)$, the chain rule implies that
	\begin{align*}
		\partial_tf(x,t)
		&=
		M^\varepsilon
		\Bigl[
		-K(t)^{-2}K'(t)\psi
		+K(t)^{-1}\partial_t\psi
		+K(t)^{-1}M^{-\varepsilon}xK'(t)\partial_x\psi
		\Bigr]
		\bigl(M^{b(t)-\varepsilon}x,t\bigr),
	\end{align*}
	and, using the expression for $K'$, we obtain
	$$
	\partial_tf(x,t)=a(t)M^{\varepsilon-b(t)}
	\psi\left(M^{b(t)-\varepsilon}x,t\right)+M^{\varepsilon-b(t)}\partial_t\psi\left(M^{b(t)-\varepsilon}x,t\right)-a(t)x\partial_x\psi\left(M^{b(t)-\varepsilon}x,t\right).
	$$
	On the other hand, $\partial_xf(x,t)=\partial_x\psi\left(M^{b(t)-\varepsilon}x,t\right)$. Consequently,
	$$
	Lf(x,t)=M^{\varepsilon-b(t)}
	\partial_t\psi\left(M^{b(t)-\varepsilon}x,t\right).
	$$
	Therefore, in order to impose the equality $Lf=h$ it suffices to choose $\psi$ so that
	$$
	\partial_t\psi(y,t)=M^{b(t)-\varepsilon}h\left(M^{\varepsilon-b(t)}y,t\right).
	$$
	Taking
	$$
	\psi(y,t)=\int_{1-N^{-\varepsilon/2}}^t
	M^{b(s)-\varepsilon}h\left(M^{\varepsilon-b(s)}y,s\right)\dd s,
	$$
	one obtains
	$$
	f(x,t)=\int_{1-N^{-\varepsilon/2}}^tM^{b(s)-b(t)}h\left(M^{b(t)-b(s)}x,s\right)\dd s,
	$$
	which satisfies $Lf=h$.
\end{proof}

Next, choose a fixed nonnegative function $\eta\in C^\infty(\mathbb{R};[0,1])$ such that $\eta(s)=0$ for $s\leq 0$ and $\eta(s)=1$ for $s\geq 1$. For $N$ sufficiently large, we have
$M^{-\varepsilon}\leq N^{-\varepsilon/2}$ by Remark~\ref{RelNM}. We may therefore define the smooth cutoff
$$
c_0(t)=\eta\bigl(M^\varepsilon(t-1+N^{-\varepsilon/2})\bigr),
\qquad t\in[0,1].
$$
By construction, $c_0$ satisfies
\begin{enumerate}
    \item $|c_0(t)|\leq 1$ for all $t\in[0,1]$,
    \item $c_0(t)=0$ for all $0\leq t\leq 1-N^{-\varepsilon/2}$ and $c_0(t)=1$ for all $1-N^{-\varepsilon/2}+M^{-\varepsilon}\leq t\leq 1,$
    \item $|c_0'(t)|\leq C M^\varepsilon$ for all $t\in[0,1]$, where $C=\|\eta'\|_{L^\infty(\mathbb{R})}$ is an absolute constant determined by the choice of $\eta$.
\end{enumerate}

\begin{rem}\label{rem:intervalofinal}
Let $t\in[1-M^{-\varepsilon/2},1]$. Since $M^{-\varepsilon}+M^{-\varepsilon/2}\leq N^{-\varepsilon/2}$ by Remark~\ref{RelNM}, we have $t\geq 1-N^{-\varepsilon/2}+M^{-\varepsilon}$, so that $c_0(t)=1$. Moreover, $0\leq\int_t^1a(s)\dd s\leq2N^{\varepsilon}M^{-\varepsilon/2}\leq\log(4/3)$ by Hypothesis~\eqref{hipotesisTHen0} and Remark~\ref{RelNM}, and therefore
\begin{equation}\label{eq:bfinal}
1\leq M^{b(t)}\leq\frac43\qquad\text{for every }t\in[1-M^{-\varepsilon/2},1].
\end{equation}
\end{rem}

Note that, as $\supp(w^1(\cdot,t))\subseteq[-1,1]$ for every $t\in[0,1]$, by Lemma~\ref{velocidadHolder} with ${R=1}$, $\|u(w^1)(\cdot,t)\|_{C^{5,1/2}}\leq C_{c(1/2)}\|w^1(\cdot,t)\|_{C^{4,1/2}}\leq C_{c(1/2)}N^6$ for every $t\in[1-N^{-\varepsilon/2},1]$, with a constant $C_{c(1/2)}$ that depends only on $c(1/2)$. We will use this observation without explicit mention again.

Motivated by Lemma~\ref{DuhamelForm}, let us consider the functions $c_1, c_2\in C^\infty([0,1])$, defined for each $t\in[0,1]$ as
	\begin{align*}
		c_1(t)&=-\int_{1-N^{-\varepsilon/2}}^t\frac{\partial^3_xu(w^1)(0,s)}{6}c_0(s)M^{2\varepsilon-2b(s)}\dd s, \\
		c_2(t)&=\int_{1-N^{-\varepsilon/2}}^{t}\partial_xw^1(0,s)c_0(s)M^{\varepsilon-b(s)}\dd s.
	\end{align*}
	Observe that $c_1(t)=c_2(t)=0$ for all $0\leq t\leq 1-N^{-\varepsilon/2}$ because $c_0(t)=0$ for all $0\leq t\leq 1-N^{-\varepsilon/2}$. Moreover, $b(t)\geq 0$ for each $t\in[1-N^{-\varepsilon/2},1]$ by Hypothesis~\eqref{hipotesisTHen0}, so that
	\begin{align*}
		|c_1(t)|&\leq \int_{1-N^{-\varepsilon/2}}^{1}\frac{\|u(w^1)(\cdot, s)\|_{C^3}}{6}|c_0(s)| M^{2\varepsilon}\dd s\leq \frac{C_{c(1/2)}}{6}N^6M^{2\varepsilon}N^{-\varepsilon/2}\leq M^{3\varepsilon}, \\
		|c_2(t)|&\leq \int_{1-N^{-\varepsilon/2}}^{1} \|w^1(\cdot, s)\|_{C^1}|c_0(s)|M^{\varepsilon}\dd s\leq N^6M^{\varepsilon}N^{-\varepsilon/2}\leq M^{2\varepsilon}.
	\end{align*}
In any case, $|c_i(t)|\leq M^{3\varepsilon}$ for each $t\in[0,1]$ and $i\in\{1,2\}$.

Recall Definition~\ref{defperfiles}. For each $i\in\{0,1,2\}$ and $(x,t)\in\R\times[0,1]$ we can define 
$$w^2_i(x,t)=c_i(t)M^{\varepsilon-b(t)}\varphi_i(M^{b(t)-\varepsilon}x), \qquad w^2(x,t)=\sum_{i=0}^{2}w^2_i(x,t).$$ Observe that, since the profiles of Definition~\ref{defperfiles} are all odd, $w^2(\cdot, t)$ is also odd for each $t\in[0,1]$.

For brevity, we will write $u^1$ instead of $u(w^1)$ and, likewise, $u^2_i$ instead of $u(w^2_i)$ for $i\in\{0,1,2\}$. Since $w^1$ already solves Equation~\eqref{eq:DGF} with forcing term $F^1$, adding the perturbation $w^2$ requires a corresponding correction to the forcing. More precisely, expanding the nonlinear terms, we find that $w^1+w^2$ solves Equation~\eqref{eq:DGF} with forcing term $F^1+F^2$, where
$$
F^2=(\partial_tw^2+u^1\partial_xw^2-w^2Hw^1)+(u^2\partial_xw^1-w^1Hw^2)+(u^2\partial_xw^2-w^2Hw^2).
$$
The terms $w^2_1$ and $w^2_2$ can be understood as corrections, which have been chosen so that, when applying Lemma~\ref{DuhamelForm}, the new force $F^2$ can be rewritten as
	\begin{align}\label{descomposicionfuerza}
		F^2&=\underbrace{c_0'(t)M^{\varepsilon-b(t)}\varphi_0(M^{b(t)-\varepsilon}x)}_{\textit{Activation term}}-\underbrace{\frac{\partial^3_xu^1(0,t)}{6}x^3\partial_xw^2_0}_{\textit{Correction given by $w^2_1$}}\\ &+\underbrace{\partial_xw^1(0,t)xc_0(t)M^{\varepsilon-b(t)}\chi(M^{b(t)-\varepsilon+1/2}x)\cos(M^{b(t)-\varepsilon+1}x)}_{\textit{Correction given by $w^2_2$}}\notag\\
		&+\underbrace{\sum_{i=0}^{2}\left[(u^1-a(t)x)\partial_xw^2_i-(Hw^1-a(t))w^2_i\right]}_{\textit{Linearized part}}+\underbrace{\sum_{i=0}^{2}\left[u^2_i\partial_xw^1-w^1Hw^2_i\right]}_{\textit{Cross terms}}\notag \\
		&+\underbrace{\sum_{i,j=0}^{2}\left[u^2_i\partial_xw^2_j-w^2_iHw^2_j\right]}_{\textit{Quadratic terms}}. \notag
	\end{align}
From the definitions of $c_0,c_1$ and $c_2$ one deduces that $w^2(\cdot,t)=F^2(\cdot, t)=0$ for all $0\leq t\leq 1-N^{-\varepsilon/2}$, that is, \eqref{eq:antesdetiempo} holds.

\subsection{Size estimates of the perturbation}\label{step2}
Observe that, by Remark~\ref{rem:perfiles} and the definition of $w^2_i$, we have
$$
\supp\bigl(w^2_i(\cdot,t)\bigr)\subseteq\{|x|\leq M^{\varepsilon-b(t)-1/2}\}
\qquad\text{for every }t\in[0,1]\text{ and }i\in\{0,1,2\}.
$$
If $t\in[1-N^{-\varepsilon/2},1]$, then $b(t)\geq0$ by Hypothesis~\eqref{hipotesisTHen0}, so this set is contained in $[-1,1]$ because $\varepsilon<1/9$. If $t\leq1-N^{-\varepsilon/2}$, then $w^2(\cdot,t)=0$. Since $\supp\bigl(w^1(\cdot,t)\bigr)\subseteq[-1,1]$ by hypothesis, we conclude that $\supp\bigl((w^1+w^2)(\cdot,t)\bigr)\subseteq[-1,1]$ for every $t\in[0,1]$. 

Observe also that, by construction of $F^2$, the same argument proves the statement about the support of the forcing term.

We next estimate the size of the profiles and of the layers.

\begin{lem}\label{tamanosol}
For every $t\in[0,1]$, $i\in\{1,2\}$ and $j\in\N_0$ the following estimates hold:
\begin{align*}
&\|\partial^j_x\varphi_0\|_{L^\infty}\lesssim_{\chi,j}M^j,
&&\|\partial^j_xw^2_0(\cdot,t)\|_{L^\infty}\lesssim_{\chi,j}M^jM^{(j-1)(b(t)-\varepsilon)},\\
&\|\partial^j_x\varphi_i\|_{L^\infty}\lesssim_{\chi,j}M^{j-1/2},
&&\|\partial^j_xw^2_i(\cdot,t)\|_{L^\infty}\lesssim_{\chi,j}M^{3\varepsilon}M^{j-1/2}M^{(j-1)(b(t)-\varepsilon)}.
\end{align*}
\end{lem}
\begin{proof}
Recall Definition~\ref{defperfiles}. In order to obtain the estimates for the profiles, apply Estimate~\eqref{eq:estima_dj_util} of Lemma~\ref{lem:estima_dj_gen}.

For the layers, the chain rule gives
$$
\partial^j_xw^2_i(x,t)=c_i(t)M^{(j-1)(b(t)-\varepsilon)}\varphi_i^{(j)}\bigl(M^{b(t)-\varepsilon}x\bigr),
$$
and the estimates follow from those for the profiles, together with $|c_0(t)|\leq1$ and $|c_i(t)|\leq M^{3\varepsilon}$ for $i\in\{1,2\}$.
\end{proof}

Estimate~\eqref{eq:reg_pert} follows readily from Lemma~\ref{tamanosol}. Indeed, let $t\in[1-M^{-\varepsilon/2},1]$. By Remark~\ref{rem:intervalofinal}, $M^{(j-1)(b(t)-\varepsilon)}\leq(4/3)^{j-1}$ for $j\geq1$ and $M^{-(b(t)-\varepsilon)}\leq M^{\varepsilon}$, so Lemma~\ref{tamanosol} gives
$$
\|\partial_x^jw^2(\cdot,t)\|_{L^\infty}\lesssim_\chi M^{3\varepsilon+j}\qquad\text{for }0\leq j\leq5.
$$
Applying Lemma~\ref{interpolacion} to $\partial_x^4w^2(\cdot,t)$ we obtain $\|w^2(\cdot,t)\|_{C^{4,1/2}}\lesssim_\chi M^{3\varepsilon+9/2}$ and therefore, by Remark~\ref{RelNM} and using that $\varepsilon<1/9$,
$$
\|(w^1+w^2)(\cdot,t)\|_{C^{4,1/2}}\leq N^6+C_\chi M^{3\varepsilon+9/2}\leq M^6,
$$
which is Estimate~\eqref{eq:reg_pert}.

\subsection{Growth of the Hilbert transform at the origin}\label{step3}
The verification of Estimate~\eqref{HIPOTESISTHEN0} will follow from the next result.
\begin{lem}\label{estimaporabajoHilbert}
There exists $C_\chi>0$ such that, for every $M > 1$ sufficiently large,
$$
\frac{1}{2}\leq H\varphi_0(0)\leq \frac{3}{2}
\qquad\text{and}\qquad
\left|H\varphi_1(0)\right|+\left|H\varphi_2(0)\right|\leq C_\chi M^{-1}.
$$
\end{lem}
\begin{proof}
By Equations~\eqref{Hphi0} and~\eqref{Hphiorigen}, using $\chi(0)=1$ we get $H\varphi_0(0)=1+r_0(0)$ and $H\varphi_i(0)=r_i(0)$ for $i\in\{1,2\}$. The bounds for $r_i$ in Corollary~\ref{lem:aglutinado} with $K=1$ and $j=0$ give $|r_0(0)|\lesssim_\chi M^{-1/2}$ and $|r_i(0)|\lesssim_\chi M^{-1}$, which proves the claim.
\end{proof}
Having established this auxiliary lemma, we can show that for every $t\in[1-M^{-\varepsilon/2},1]$,
$$
\frac{M^\varepsilon}{4}\leq H(w^1+w^2)(0,t)\leq 2M^\varepsilon.
$$
Indeed, fix such a $t\in[1-M^{-\varepsilon/2},1]$. By Remark~\ref{rem:intervalofinal}, $c_0(t)=1$ and $\frac34M^{\varepsilon}\leq M^{\varepsilon-b(t)}\leq M^{\varepsilon}$. Since the Hilbert transform commutes with dilations, $Hw^2_i(0,t)=c_i(t)M^{\varepsilon-b(t)}H\varphi_i(0)$, and therefore
$$
H(w^1+w^2)(0,t)=a(t)+M^{\varepsilon-b(t)}H\varphi_0(0)+\sum_{i=1}^2c_i(t)M^{\varepsilon-b(t)}H\varphi_i(0).
$$
By Lemma~\ref{estimaporabajoHilbert}, which applies since $M\geq N$ is large, the last sum is bounded in absolute value by $C_\chi M^{4\varepsilon-1}\leq M^{\varepsilon}/8$, where the last inequality follows because $\varepsilon<1/9$. Using Hypothesis~\eqref{hipotesisTHen0}, $\frac12\leq H\varphi_0(0)\leq\frac32$ and $2N^{\varepsilon}\leq M^{\varepsilon}/8$, we conclude that
$$
\frac{M^\varepsilon}{4}=\frac38M^{\varepsilon}-\frac18M^{\varepsilon}\leq H(w^1+w^2)(0,t)\leq\frac18M^{\varepsilon}+\frac32M^{\varepsilon}+\frac18M^{\varepsilon}\leq2M^{\varepsilon},
$$
which is Estimate~\eqref{HIPOTESISTHEN0}.

\subsection{Estimates of the forcing term}\label{step4}
We will establish Inequality~\eqref{eq:estimafuerza} via the following series of results.
\subsubsection{Bounds for the activation term}
\begin{lem}\label{acotacionSusto}
Let $0<\alpha<1$ and $t\in[0,1]$. Then
$$
\|c_0'(t)M^{\varepsilon-b(t)}\varphi_0(M^{b(t)-\varepsilon}\cdot)\|_{C^{0,\alpha}}\lesssim_{\chi,\varepsilon} M^{2\varepsilon+(1+1/\varepsilon)\alpha-1/\varepsilon}.
$$
\end{lem}

\begin{proof}
By assumption, $N$ is sufficiently large and, by definition,
$$
|c_0'(t)|\leq CM^\varepsilon\ind_{[1-N^{-\varepsilon/2},1-N^{-\varepsilon/2}+M^{-\varepsilon}]}(t).
$$
Moreover,
$$
\partial_x[c_0'(t)M^{\varepsilon-b(t)}\varphi_0(M^{b(t)-\varepsilon}x)]=c_0'(t)\varphi_0'(M^{b(t)-\varepsilon}x),
$$
and hence
$$
\|\partial_x[c_0'(t)M^{\varepsilon-b(t)}\varphi_0(M^{b(t)-\varepsilon}\cdot)]\|_{L^\infty}\lesssim_\chi M^{\varepsilon+1}.
$$
Since $a(\tau)>0$ for every $\tau\geq 1-N^{-\varepsilon/2}$, the function $K$ is decreasing on the interval $[1-N^{-\varepsilon/2},1-N^{-\varepsilon/2}+M^{-\varepsilon}]$. Thus, since $c_0'=0$ outside this interval and $M^{b(\tau)}=K(\tau)$, we obtain
$$
\begin{aligned}
\|c_0'(t)M^{\varepsilon-b(t)}\varphi_0(M^{b(t)-\varepsilon}\cdot)\|_{L^\infty}
&\leq \sup_{\tau\in[1-N^{-\varepsilon/2},1-N^{-\varepsilon/2}+M^{-\varepsilon}]}CM^\varepsilon M^{\varepsilon-b(\tau)} \\
&\leq \frac{CM^{2\varepsilon}}{K(1-N^{-\varepsilon/2}+M^{-\varepsilon})}.
\end{aligned}
$$
On the other hand, since $a(\tau)\leq 2N^\varepsilon$ for every $\tau\geq 1-N^{-\varepsilon/2}$, Remark~\ref{RelNM} yields
\begin{align*}
M^{1/\varepsilon}
&=e^{\int_{1-N^{-\varepsilon/2}}^1a(s)\dd s}=e^{\int_{1-N^{-\varepsilon/2}}^{1-N^{-\varepsilon/2}+M^{-\varepsilon}}a(s)\dd s}K(1-N^{-\varepsilon/2}+M^{-\varepsilon}) \\
&\leq e^{2N^\varepsilon M^{-\varepsilon}}K(1-N^{-\varepsilon/2}+M^{-\varepsilon})\leq 2K(1-N^{-\varepsilon/2}+M^{-\varepsilon}).
\end{align*}
Consequently,
$$
\|c_0'(t)M^{\varepsilon-b(t)}\varphi_0(M^{b(t)-\varepsilon}\cdot)\|_{L^\infty}\leq 2CM^{2\varepsilon-1/\varepsilon}.
$$
We can now apply Lemma~\ref{interpolacion} to obtain
$$
\|c_0'(t)M^{\varepsilon-b(t)}\varphi_0(M^{b(t)-\varepsilon}\cdot)\|_{C^{0,\alpha}}\lesssim_{\chi,\varepsilon} M^{2\varepsilon-1/\varepsilon}(1+M^{\alpha/\varepsilon-2\varepsilon\alpha}M^{\varepsilon\alpha+\alpha})\lesssim_\chi M^{2\varepsilon+(1+1/\varepsilon)\alpha-1/\varepsilon}.
$$
\end{proof}

\subsubsection{Bounds for the quadratic term}
Among all the terms in the force decomposition~\eqref{descomposicionfuerza}, the most challenging one to estimate is the quadratic term, since we essentially need to bound nine different terms. Recall that the quadratic term is given by
$$
\begin{aligned}
\sum_{i,j=0}^2[u(w^2_i)\partial_xw^2_j-w^2_iHw^2_j](x,t)
&=\sum_{i,j=0}^2c_i(t)M^{\varepsilon-b(t)}
\left(\int_0^xH\varphi_i(M^{b(t)-\varepsilon}y)\dd y\right)
c_j(t)\varphi'_j(M^{b(t)-\varepsilon}x)\\
&-\sum_{i,j=0}^2c_i(t)c_j(t)M^{2\varepsilon-2b(t)}
\varphi_i(M^{b(t)-\varepsilon}x)H\varphi_j(M^{b(t)-\varepsilon}x).
\end{aligned}
$$
By a change of variables, we can rewrite the quadratic term as follows:
\begin{equation*}
\sum_{i,j=0}^2[u(w^2_i)\partial_xw^2_j-w^2_iHw^2_j](x,t)
=\sum_{i,j=0}^2c_i(t)c_j(t)M^{2\varepsilon-2b(t)}f_{ij}(M^{b(t)-\varepsilon}x),
\end{equation*}
where, for each $i,j\in\{0,1,2\}$,
$$
f_{ij}(z)=\left(\int_0^zH\varphi_i(y)\dd y\right)\varphi_j'(z)-\varphi_i(z)H\varphi_j(z).
$$
Although one could attempt to estimate these nine functions directly, we shall employ Equations~\eqref{reescrituraperfiles} and~\eqref{reescrituraperfilesP} from Corollary~\ref{lem:aglutinado} in order to make explicit the key cancellation arising from the specific functional form of the profiles introduced in Definition~\ref{defperfiles}.

\begin{rem}[Fundamental cancellation]\label{cancelacionfund}
By using Equations~\eqref{reescrituraperfiles} and~\eqref{reescrituraperfilesP} from Corollary~\ref{lem:aglutinado}, observe that, for any $i,j\in\{0,1,2\}$,
\begin{align*}
f_{ij}&=\left(-\frac{\varphi_i}{M}+E_i+s_i\right)\varphi_j'-\varphi_i\left(\frac{-\varphi_j'}{M}+e_j+r_j\right)\\
&=\cancel{-\frac{\varphi_i\varphi_j'}{M}}+E_i\varphi_j'+s_i\varphi_j'+\cancel{\frac{\varphi_i\varphi_j'}{M}}-\varphi_ie_j-\varphi_ir_j=E_i\varphi_j'+s_i\varphi_j'-\varphi_ie_j-\varphi_ir_j,
\\
\partial_xf_{ij}&=E_i'\varphi_j'+E_i\varphi_j''+s_i'\varphi_j'+s_i\varphi_j''-\varphi_i'e_j-\varphi_ie_j'-\varphi_i'r_j-\varphi_ir_j'.
\end{align*}
The terms that cancel are larger by a factor of $M^{1/2}$ than the remaining terms. This cancellation is the key mechanism underlying the proofs of the following lemmas.
\end{rem}

\begin{lem}[Bound for the quadratic term I]\label{acotacionCuad1}
		Let $i,j\in\{1,2\}$. Then
		$$
		\|f_{ij}\|_{L^\infty}\lesssim_{\chi}M^{-3/2}, \qquad \|f_{i0}\|_{L^\infty}\lesssim_\chi M^{-1}, \qquad \|f_{0j}\|_{L^\infty}\lesssim_\chi M^{-1},
		$$
		$$
     \|\partial_xf_{ij}\|_{L^\infty}\lesssim_{\chi}M^{-1/2},\qquad \|\partial_xf_{i0}\|_{L^\infty}\lesssim_{\chi}1,\qquad  \|\partial_xf_{0j}\|_{L^\infty}\lesssim_{\chi}1.
		$$
\end{lem}
\begin{proof}
Fix $i,j\in\{1,2\}$. By Remark~\ref{cancelacionfund},
\begin{equation}\label{eq:fij}
    f_{ij}=E_i\varphi_j'+s_i\varphi_j'-\varphi_ie_j-\varphi_ir_j.
\end{equation}
By Lemma~\ref{tamanosol}, $\|\varphi_i\|_{L^\infty}\lesssim_\chi M^{-1/2}$ and $\|\varphi_j'\|_{L^\infty}\lesssim_\chi M^{1/2}$, while Corollary~\ref{lem:aglutinado} yields $\|E_i\|_{L^\infty}\lesssim_\chi M^{-2}$, $\|e_j\|_{L^\infty}\lesssim_\chi M^{-1}$, $\|r_j\|_{L^\infty}\lesssim_\chi M^{-3/2}$ and, taking $L=M^{-1/2}$, we also have $\|s_i\|_{L^\infty(\{|x|\leq M^{-1/2}\})}\lesssim_\chi M^{-2}$. These estimates give the desired bound.

Indeed, the latter estimate applies since $\supp(\varphi_j')\subseteq\{|x|\leq M^{-1/2}\}$. The four terms in the expression for $f_{ij}$ are, respectively, of order $\mathcal{O}_\chi(M^{-3/2})$, $\mathcal{O}_\chi(M^{-3/2})$, $\mathcal{O}_\chi(M^{-3/2})$, and $\mathcal{O}_\chi(M^{-2})$. Consequently, $\|f_{ij}\|_{L^\infty}\lesssim_\chi M^{-3/2}.$
Differentiating in Equation~\eqref{eq:fij}, we obtain
\begin{equation}\label{eq:dfij}
    \partial_xf_{ij}
=E_i'\varphi_j'+E_i\varphi_j''
+s_i'\varphi_j'+s_i\varphi_j''-\varphi_i'e_j-\varphi_ie_j'
-\varphi_i'r_j-\varphi_ir_j'.
\end{equation}
Proceeding in the same way, and using the additional estimates
$$
\|\varphi_i'\|_{L^\infty}\lesssim_\chi M^{1/2},
\qquad
\|\varphi_j''\|_{L^\infty}\lesssim_\chi M^{3/2},
\qquad
\|E_i'\|_{L^\infty}\lesssim_\chi M^{-1},
$$
$$
\|s_i'\|_{L^\infty(\{|x|\leq M^{-1/2}\})}\lesssim_\chi M^{-3/2},
\qquad
\|e_j'\|_{L^\infty}\lesssim_\chi 1,
$$
as well as the corresponding estimates for the remainder terms $r_j$ and $r_j'$, we find that each of the eight terms is of order $\mathcal{O}_\chi(M^{-1/2})$. Therefore, $\|\partial_xf_{ij}\|_{L^\infty}\lesssim_\chi M^{-1/2}.$

Finally, the cases $f_{i0}$ and $f_{0j}$ are handled in exactly the same manner, using the corresponding estimates for $\varphi_0$, $e_0$, $E_0$, and $r_0$ from Corollary~\ref{lem:aglutinado} and Lemma~\ref{tamanosol}.
\end{proof}
The function $f_{00}$ is somewhat special. The approach used in Lemma~\ref{acotacionCuad1} is too crude, as it fails to exploit a certain ``frequency cancellation''. This is the subject of the following result. Observe that
\begin{align*}
f_{00}(x)&=\left(-\frac{\varphi_0(x)}{M}+\frac{\chi'(M^{1/2}x)\cos(Mx)}{M^{3/2}}-\frac{\chi''(M^{1/2}x)\sin(Mx)}{M^2}+s_0(x)\right)\varphi_0'(x)\\
&-\varphi_0(x)\left(\chi(M^{1/2}x)\cos(Mx)+r_0(x)\right)\\
&=(s_0\varphi_0'-\varphi_0r_0)(x)
-\underbrace{\frac{\chi(M^{1/2}x)\chi'(M^{1/2}x)}{M^{1/2}}}_{I}
-\underbrace{\frac{\left(\chi'(M^{1/2}x)\right)^2\sin(Mx)\cos(Mx)}{M}}_{II}\\
&+\underbrace{\frac{\chi'(M^{1/2}x)\chi''(M^{1/2}x)\sin^2(Mx)}{M^{3/2}}}_{III}
+\underbrace{\frac{\chi(M^{1/2}x)\chi''(M^{1/2}x)\sin(Mx)\cos(Mx)}{M}}_{IV}.
\end{align*}
\begin{lem}[Bounds for the quadratic term II]\label{acotacionCuad2}
The following estimates hold.
\begin{align*}
\|I\|_{L^\infty}&\lesssim_\chi M^{-1/2}, & \|\partial_xI\|_{L^\infty}&\lesssim_\chi 1,\\
\|II\|_{L^\infty}&\lesssim_\chi M^{-1}, & \|\partial_xII\|_{L^\infty}&\lesssim_\chi 1,\\
\|III\|_{L^\infty}&\lesssim_\chi M^{-3/2}, & \|\partial_xIII\|_{L^\infty}&\lesssim_\chi M^{-1/2},\\
\|IV\|_{L^\infty}&\lesssim_\chi M^{-1}, & \|\partial_xIV\|_{L^\infty}&\lesssim_\chi 1,\\
\|s_0\varphi_0'\|_{L^\infty}&\lesssim_\chi M^{-1}, & \|\partial_x(s_0\varphi_0')\|_{L^\infty}&\lesssim_\chi 1,\\
\|\varphi_0r_0\|_{L^\infty}&\lesssim_\chi M^{-3/2}, & \|\partial_x(\varphi_0r_0)\|_{L^\infty}&\lesssim_\chi M^{-1/2}.
\end{align*}
\end{lem}
\begin{proof}
The first four lines of estimates are elementary. The remaining ones follow by applying Estimates~\eqref{cotasr} and~\eqref{cotass} with $L=M^{-1/2}$ together with Lemma~\ref{tamanosol}.
\end{proof}
\begin{lem}[Bounds for the quadratic term III]\label{acotacionCuad3}
Let $0<\alpha<1$ and $t\in[1-N^{-\varepsilon/2},1]$. Then
$$
\|u^2\partial_xw^2-w^2Hw^2\|_{C^{0,\alpha}}\lesssim_{\chi}M^{8\varepsilon+\alpha-3/2}+M^{5\varepsilon+\alpha-1}+M^{2\varepsilon+\frac{\alpha-1}{2}}.
$$
\end{lem}

\begin{proof}
Fix $i,j\in\{0,1,2\}$. By Lemma~\ref{interpolacion} we have
\begin{equation}\label{desfinalcuad}
				\|u^2_i\partial_xw^2_j-w^2_iHw^2_j\|_{C^{0,\alpha}}\leq M^{2\varepsilon-2b(t)}|c_i(t)c_j(t)|(\|f_{ij}\|_{L^\infty}+2^{1-\alpha}M^{\alpha(b(t)-\varepsilon)}\|f_{ij}\|_{L^\infty}^{1-\alpha}\|\partial_xf_{ij}\|_{L^\infty}^\alpha).
			\end{equation}
Note that $1<2-\alpha<2$ and $b(t)\geq0$ for every $t\in[1-N^{-\varepsilon/2},1]$. Consequently, $M^{(2-\alpha)(\varepsilon-b(t))}\leq M^{2\varepsilon}.$
The proof is now completed by applying Lemmas~\ref{acotacionCuad1} and~\ref{acotacionCuad2} to the corresponding cases in estimate~\eqref{desfinalcuad}, and combining the resulting bounds. This yields
$$
\|u^2\partial_xw^2-w^2Hw^2\|_{C^{0,\alpha}}\lesssim_{\chi}M^{8\varepsilon+\alpha-3/2}+M^{5\varepsilon+\alpha-1}+M^{2\varepsilon+\frac{\alpha-1}{2}},
$$
as claimed.
\end{proof}

\subsubsection{Bounds for the linearized term}
We will use the following lemma to exploit two key features of our construction, namely, that the perturbation is much more localized in space than the background layer, and that cancellations due to parity lead to improved estimates.

\begin{lem}\label{lem:producto}
Let $p\in\N$, $0<\rho\leq1$, $h\in C^p([-1,1])$ with $h^{(q)}(0)=0$ for every $0\leq q<p$, and $g\in C^1(\R)$ with $\supp(g)\subseteq[-\rho,\rho]$. Then
\begin{align*}
\|hg\|_{L^\infty}&\leq\rho^p\|h^{(p)}\|_{L^\infty([-1,1])}\|g\|_{L^\infty},\\
\|\partial_x(hg)\|_{L^\infty}&\leq\rho^{p-1}\|h^{(p)}\|_{L^\infty([-1,1])}\bigl(\|g\|_{L^\infty}+\rho\|g'\|_{L^\infty}\bigr).
\end{align*}
\end{lem}
\begin{proof}
By Taylor's theorem, $|h(x)|\leq\|h^{(p)}\|_{L^\infty([-1,1])}|x|^p$ and $|h'(x)|\leq\|h^{(p)}\|_{L^\infty([-1,1])}|x|^{p-1}$ for each $x\in[-1,1]$. Both estimates follow from this because $g$ and $g'$ vanish outside $[-\rho,\rho]$.
\end{proof}

\begin{lem}[Bounds for the linearized term]\label{acotacionLin}
Let $0<\alpha<1$ and $t\in[1-N^{-\varepsilon/2},1]$. Then the following estimates hold.
\begin{align}
    \label{eqacotacionlin1}
\|(u^1-a(t)x)\partial_xw^2_i\|_{C^{0,\alpha}}&\lesssim_{\chi,\varepsilon} M^{7\varepsilon+\alpha-1}, \quad i\in\{1,2\}, \\
\label{eqacotacionlin2}
\left\|\left(u^1-a(t)x-\frac{\partial^3_xu^1(0,t)}{6}x^3\right)\partial_xw^2_0\right\|_{C^{0,\alpha}}&\lesssim_{\chi,\varepsilon} M^{6\varepsilon+\alpha-3/2}, \\
\label{eqacotacionlin3}
\|(Hw^1-a(t))w^2_0\|_{C^{0,\alpha}}&\lesssim_{\chi,\varepsilon} M^{4\varepsilon+\alpha-1}, \\
\label{eqacotacionlin4}
\|(Hw^1-a(t))w^2_i\|_{C^{0,\alpha}}&\lesssim_{\chi,\varepsilon} M^{7\varepsilon+\alpha-3/2}, \quad i\in\{1,2\}.
\end{align}
\end{lem}
\begin{proof}
Fix $t\in[1-N^{-\varepsilon/2},1]$, so that $b(t)\geq0$. Since $w^1(\cdot,t)$ is odd, $Hw^1(\cdot,t)$ is even and $u^1(\cdot,t)$ is odd. Hence the even-order derivatives of $u^1(\cdot,t)$ and the odd-order derivatives of $Hw^1(\cdot,t)$ vanish at the origin, while $\partial_xu^1(0,t)=Hw^1(0,t)=a(t)$. Consequently, the functions
$$
x\mapsto u^1(x,t)-a(t)x,\qquad x\mapsto u^1(x,t)-a(t)x-\frac{\partial_x^3u^1(0,t)}{6}x^3,\qquad x\mapsto Hw^1(x,t)-a(t),
$$
satisfy the hypotheses on $h$ in Lemma~\ref{lem:producto} with $p=3$, $p=5$ and $p=2$, respectively. Their $p$-th derivatives are $\partial_x^3u^1(\cdot,t)$, $\partial_x^5u^1(\cdot,t)$ and $\partial_x^3u^1(\cdot,t)$, whose norms in $L^\infty([-1,1])$ are bounded by $\|u^1(\cdot,t)\|_{C^5}\lesssim N^6\leq M^\varepsilon$, by Lemma~\ref{velocidadHolder} and Remark~\ref{RelNM}.

On the other hand, $\supp(w^2_i(\cdot,t))\subseteq\{|x|\leq M^{\varepsilon-b(t)-1/2}\}\subseteq[-1,1]$ for every $i\in\{0,1,2\}$, as shown in Subsection~\ref{step2}. We may therefore apply Lemma~\ref{lem:producto} with $g=w^2_i(\cdot,t)$ or $g=\partial_xw^2_i(\cdot,t)$ and $\rho=M^{\varepsilon-b(t)-1/2}\leq 1$.

Combined with Lemma~\ref{tamanosol} and $b(t)\geq0$, this shows that, up to constants depending only on $\chi$ and $\varepsilon$, the $L^\infty$ norms of the functions in~\eqref{eqacotacionlin1}--\eqref{eqacotacionlin4} are bounded by $M^{7\varepsilon-1}$, $M^{6\varepsilon-3/2}$, $M^{4\varepsilon-1}$ and $M^{7\varepsilon-3/2}$, respectively, and those of their derivatives by $M^{1-\varepsilon}$ times the same quantities. Lemma~\ref{interpolacion} then yields Estimates~\eqref{eqacotacionlin1}--\eqref{eqacotacionlin4}.
\end{proof}

\subsubsection{Bounds for the cross terms}
We now estimate the remaining cross terms.

\begin{lem}[Bounds for the cross terms]\label{acotacionCruzI}
Let $0<\alpha<1$ and $t\in[1-N^{-\varepsilon/2},1]$. Then the following estimates hold.
\begin{align}
    \label{eq:acotacioncruz1}
\|w^1Hw^2_i\|_{C^{0,\alpha}}&\lesssim_{\chi,\varepsilon} M^{6\varepsilon+\alpha-1}, \quad i\in\{1,2\}, \\
\label{eq:acotacioncruz2}
\left\|-w^1Hw^2_0+\partial_xw^1(0,t)x c_0(t)M^{\varepsilon-b(t)}
\Phi(M^{b(t)-\varepsilon}x)\right\|_{C^{0,\alpha}}
&\lesssim_{\chi,\varepsilon} M^{5\varepsilon+\alpha-3/2},\\
\label{eq:acotacioncruz3}
\|u^2_0\partial_xw^1\|_{C^{0,\alpha}}&\lesssim_{\chi,\varepsilon}M^{3\varepsilon+\alpha-1}, \\
\label{eq:acotacioncruz4}
\|u^2_i\partial_xw^1\|_{C^{0,\alpha}}&\lesssim_{\chi,\varepsilon}M^{6\varepsilon+\alpha-3/2}, \quad i\in\{1,2\},
\end{align}
where $\Phi(x)=\chi(M^{1/2}x)\cos(Mx)$.
\end{lem}

\begin{proof}
Fix $t\in[1-N^{-\varepsilon/2},1]$, so that $b(t)\geq0$. Recall that $w^1(\cdot,t)$ is odd, $\supp(w^1(\cdot,t))\subseteq[-1,1]$ and $\|w^1(\cdot,t)\|_{C^{4,1/2}}\leq N^6\leq M^\varepsilon$ by Remark~\ref{RelNM}, and that $|c_0(t)|\leq1$ and $|c_i(t)|\leq M^{3\varepsilon}$ for $i\in\{1,2\}$.

Let $i\in\{0,1,2\}$. Since the Hilbert transform commutes with dilations, $Hw^2_i(x,t)=c_i(t)M^{\varepsilon-b(t)}(H\varphi_i)(M^{b(t)-\varepsilon}x)$, and Corollary~\ref{lem:aglutinado} yields
\begin{equation}\label{eq:reescrituracruz1}
    Hw^2_i(x,t)=c_i(t)M^{\varepsilon-b(t)}\left(-\frac{\varphi_i'}{M}+e_i\right)(M^{b(t)-\varepsilon}x)+c_i(t)M^{\varepsilon-b(t)}r_i(M^{b(t)-\varepsilon}x).
\end{equation}

We begin by proving Estimate~\eqref{eq:acotacioncruz1}. Fix $i\in\{1,2\}$. We study each of the summands in Equation~\eqref{eq:reescrituracruz1} separately. The one involving $r_i$ is easily handled: Estimate~\eqref{cotasr} with $K=2$ gives $\|r_i\|_{L^\infty}\lesssim_\chi M^{-3/2}$ and $\|\partial_xr_i\|_{L^\infty}\lesssim_\chi M^{-1}$, and the corresponding bound follows from Lemma~\ref{interpolacion}. As for the other one, observe that, by Corollary~\ref{lem:aglutinado} and Lemma~\ref{tamanosol},
$$
\left\|-\frac{\varphi_i'}{M}+e_i\right\|_{L^\infty}\lesssim_{\chi} M^{-1/2}, \quad \left\|\partial_x\left(-\frac{\varphi_i'}{M}+e_i\right)\right\|_{L^\infty}\lesssim_{\chi} M^{1/2}.
$$
Moreover, the first summand in Equation~\eqref{eq:reescrituracruz1} is supported on $\{|x|\leq M^{\varepsilon-b(t)-1/2}\}$ by Remark~\ref{rem:perfiles} and Corollary~\ref{lem:aglutinado}. Therefore, by applying Lemma~\ref{lem:producto} with $h=w^1(\cdot, t)$, $p=1$, $\rho=M^{\varepsilon-b(t)-1/2}$ and $g=c_i(t)M^{\varepsilon-b(t)}\left(-\frac{\varphi_i'}{M}+e_i\right)(M^{b(t)-\varepsilon}\cdot)$ together with Lemma~\ref{interpolacion}, Estimate~\eqref{eq:acotacioncruz1} follows.

We next prove Estimate~\eqref{eq:acotacioncruz2}. Comparing the two expressions for $H\varphi_0$ provided by Corollary~\ref{lem:aglutinado}, we see that $-\varphi_0'/M+e_0=\Phi$. Hence, by Equation~\eqref{eq:reescrituracruz1} with $i=0$, the function in~\eqref{eq:acotacioncruz2} equals
$$
-\bigl(w^1(x,t)-\partial_xw^1(0,t)x\bigr)c_0(t)M^{\varepsilon-b(t)}\Phi(M^{b(t)-\varepsilon}x)-w^1(x,t)c_0(t)M^{\varepsilon-b(t)}r_0(M^{b(t)-\varepsilon}x).
$$
The second term is handled as before, now using Estimate~\eqref{cotasr} with $K=3$, which gives $\|r_0\|_{L^\infty}\lesssim_\chi M^{-3/2}$ and $\|\partial_xr_0\|_{L^\infty}\lesssim_\chi M^{-1}$. As for the first one, since $w^1(\cdot,t)$ is odd, the function $w^1(\cdot,t)-\partial_xw^1(0,t)(\cdot)$ vanishes at the origin together with its first two derivatives, and its third derivative is $\partial_x^3w^1(\cdot,t)$. Moreover, as $\Phi=\operatorname{Re}\left(g^{(0,0)}_{M^{1/2}, M}\right)$, Lemma~\ref{lem:estima_dj_gen} gives $\|\partial^j_x\Phi\|_{L^\infty}\lesssim_\chi M^j$ and, by definition, $\Phi$ is supported on $\{|x|\leq M^{-1/2}\}$. Therefore, by applying Lemma~\ref{lem:producto} with $h=w^1(\cdot,t)-\partial_xw^1(0,t)(\cdot)$, $p=3$, $\rho=M^{\varepsilon-b(t)-1/2}$ and $g=c_0(t)M^{\varepsilon-b(t)}\Phi(M^{b(t)-\varepsilon}\cdot)$ together with Lemma~\ref{interpolacion}, Estimate~\eqref{eq:acotacioncruz2} follows.

Finally, we prove Estimates~\eqref{eq:acotacioncruz3} and~\eqref{eq:acotacioncruz4}. Let $i\in\{0,1,2\}$. By a change of variables,
$$u^2_i(x,t)=c_i(t)M^{2\varepsilon-2b(t)}\int_0^{M^{b(t)-\varepsilon}x}H\varphi_i(y)\dd y.$$
Corollary~\ref{lem:aglutinado} yields
\begin{equation}\label{eq:reescrituracruz2}
    u^2_i(x,t)=c_i(t)M^{2\varepsilon-2b(t)}\left(-\frac{\varphi_i}{M}+E_i\right)(M^{b(t)-\varepsilon}x)+c_i(t)M^{2\varepsilon-2b(t)}s_i(M^{b(t)-\varepsilon}x).
\end{equation}
Again, we study each of the summands separately. The one involving $s_i$ is handled as before. Since $\partial_xw^1(\cdot,t)$ vanishes outside $[-1,1]$, only the values of $s_i$ on $\{|x|\leq M^{b(t)-\varepsilon}\}$ are involved, and Estimate~\eqref{cotass} with $L=M^{b(t)-\varepsilon}$ gives $\|s_i\|_{L^\infty(\{|x|\leq M^{b(t)-\varepsilon}\})}\lesssim_\chi M^{b(t)-\varepsilon-3/2}$ and $\|\partial_xs_i\|_{L^\infty}\lesssim_\chi M^{-3/2}$, so the corresponding bound follows from Lemma~\ref{interpolacion}. As for the other one, no cancellation is needed, since taking the primitive already yields a factor $M^{-1}$. Indeed, by Corollary~\ref{lem:aglutinado} and Lemma~\ref{tamanosol}, for $j\in\{0,1\}$ we have
$$
\left\|\partial_x^j\left(-\frac{\varphi_0}{M}+E_0\right)\right\|_{L^\infty}\lesssim_\chi M^{j-1},\qquad\left\|\partial_x^j\left(-\frac{\varphi_i}{M}+E_i\right)\right\|_{L^\infty}\lesssim_\chi M^{j-3/2}\quad\text{for }i\in\{1,2\},
$$
and Estimates~\eqref{eq:acotacioncruz3} and~\eqref{eq:acotacioncruz4} follow from the Leibniz rule and Lemma~\ref{interpolacion}.
\end{proof}

By Property~\eqref{eq:antesdetiempo}, the desired estimate for the forcing term is immediate on the time interval $[0,1-N^{-\varepsilon/2}]$. It therefore remains to consider $t\in[1-N^{-\varepsilon/2},1]$. Choosing $\alpha=1-9\varepsilon$, the assumption $0<\varepsilon<1/9$ together with Lemmas~\ref{acotacionSusto}, \ref{acotacionCuad3}, \ref{acotacionLin} and \ref{acotacionCruzI} yields
$$
\|F^2(\cdot,t)\|_{C^{0,1-9\varepsilon}}\leq C_{\chi,\varepsilon}M^{-2\varepsilon},
$$
where $C_{\chi,\varepsilon}>0$ depends only on $c(1/2)$, $\|\chi\|_{C^5}$, $\varepsilon$, and the constant appearing in Lemma~\ref{auxiliarwave}. Combining both time intervals, we conclude that
\begin{equation*}
\sup_{0\leq t\leq1}\|F^2(\cdot,t)\|_{C^{0,1-9\varepsilon}}
\leq C_{\chi,\varepsilon}M^{-2\varepsilon},
\end{equation*}
that is, \eqref{eq:estimafuerza} holds. That $F^2\in C([0,1]; C^{0,1-9\varepsilon}(\R))$ actually follows from the fact that $F^2\in C^\infty(\R\times[0,1])$ and $\supp(F^2(\cdot, t))\subseteq[-1,1]$ for every $t\in[0,1]$. This finishes the proof of Theorem~\ref{thmIterativo}.
\end{proof}

\section{Constructing the solution}

The statement of Theorem~\ref{thmIterativo} already suggests how the desired solution to Equation~\eqref{eq:DGF} will be constructed: we begin with a smooth, compactly supported initial layer satisfying the prescribed hypotheses, and then proceed inductively by adding successive layers at appropriately chosen times.

Although the computations are similar to those carried out in parts of the previous section, for the sake of completeness we first check that the layer from which the iteration will start satisfies the hypotheses of Theorem~\ref{thmIterativo}.
\begin{lem}\label{arrancar_ind}
Let $0<\varepsilon<1/9$ and $N>1$ sufficiently large. Then the functions defined for $(x,t)\in\R\times[0,1]$ by
$$
w^{(0)}(x,t)=-N^\varepsilon\chi(N^{1/2-\varepsilon}x)\sin(N^{1-\varepsilon}x)
$$
and
$$
F^{(0)}(x,t)=\left[u(w^{(0)})\partial_xw^{(0)}-w^{(0)}Hw^{(0)}\right](x,t)
$$
have the following properties:
\begin{enumerate}
    \item\label{capa0:sol} $w^{(0)}$ and $F^{(0)}$ belong to $C^\infty(\R\times[0,1])$, do not depend on $t$, and $w^{(0)}$ solves Equation~\eqref{eq:DGF} with forcing term $F^{(0)}$.
    \item\label{capa0:sop} $w^{(0)}(\cdot,t)$ is odd and
    $$
    \supp\bigl(w^{(0)}(\cdot,t)\bigr)\cup\supp\bigl(F^{(0)}(\cdot,t)\bigr)\subseteq[-N^{\varepsilon-1/2},N^{\varepsilon-1/2}]\subseteq[-1,1]
    $$
    for every $t\in[0,1]$.
    \item\label{capa0:hol} $F^{(0)}\in C\bigl([0,1];C^{0,\alpha}(\R)\bigr)$ for every $\alpha\in(0,1)$.
    \item\label{capa0:est} For every $t\in[0,1]$,
    $$
    \|w^{(0)}(\cdot,t)\|_{C^{4,1/2}}\leq N^6
    \qquad\text{and}\qquad
    \frac{N^\varepsilon}{2}\leq Hw^{(0)}(0,t)\leq\frac{3N^\varepsilon}{2}.
    $$
\end{enumerate}
In particular, $w^{(0)}$ and $F^{(0)}$ satisfy all the hypotheses imposed on $w^1$ and $F^1$ in Theorem~\ref{thmIterativo} for this value of $N$.
\end{lem}
\begin{proof}
Since neither $w^{(0)}$ nor $F^{(0)}$ depends on $t$, we omit the time variable throughout the proof.
 
We begin by proving Items~\ref{capa0:sol} and~\ref{capa0:sop}. Note that $w^{(0)}$ is odd and $w^{(0)}\in C_c^\infty(\R)$. Recall that $\supp(\chi)\subseteq[-1,1]$, from which it follows that $\supp(w^{(0)})\subseteq[-N^{\varepsilon-1/2},N^{\varepsilon-1/2}]\subseteq[-1,1].$

As $w^{(0)}\in C^\infty_c(\R)$ and the Hilbert transform is bounded on $H^s(\R)$ for every $s\geq0$, we have $Hw^{(0)}\in C^\infty(\R)$, and therefore its primitive $u(w^{(0)})$ is smooth as well. Consequently $F^{(0)}$ is infinitely differentiable, and since each of its two terms contains either the factor $w^{(0)}$ or its derivative, we get $\supp(F^{(0)})\subseteq\supp(w^{(0)})$. 

Finally, observe that since $\partial_tw^{(0)}=0$, $F^{(0)}$ has been chosen so that the identity
$$
\partial_tw^{(0)}+u(w^{(0)})\partial_xw^{(0)}=w^{(0)}Hw^{(0)}+F^{(0)}
$$
holds, that is, $w^{(0)}$ solves Equation~\eqref{eq:DGF} with forcing term $F^{(0)}$.
 
About Item~\ref{capa0:hol}, note that by the previous step $F^{(0)}\in C_c^\infty(\R)\subseteq C^1(\R)$, so Lemma~\ref{interpolacion} gives $F^{(0)}\in C^{0,\alpha}(\R)$ for every $\alpha\in(0,1)$. Since $F^{(0)}$ does not depend on $t$, the map $t\mapsto F^{(0)}(\cdot,t)$ is constant and, in particular, continuous from $[0,1]$ into $C^{0,\alpha}(\R)$.
 
In order to prove Item~\ref{capa0:est} we begin with the bound for the $C^{4,1/2}$ norm. Let $0\leq j\leq5$. Observe that $w^{(0)}=-N^\varepsilon\operatorname{Im}\left(g^{(0,0)}_{N^{1/2-\varepsilon}, N^{1-\varepsilon}}\right)$, so that Estimate~\eqref{eq:estima_dj_util} of Lemma~\ref{lem:estima_dj_gen} gives that
$$
\|\partial^j_x w^{(0)}\|_{L^\infty}\leq C_\chi N^\varepsilon N^{j(1-\varepsilon)}\leq C_\chi N^5,
$$
for some $C_\chi>0$. Therefore $\|w^{(0)}\|_{C^4}\leq 5C_\chi N^5$ and, by Lemma~\ref{interpolacion} with $\alpha=1/2$,
$$
[\partial_x^4w^{(0)}]_{C^{0,1/2}}
\leq\sqrt2\|\partial_x^4w^{(0)}\|_{L^\infty}^{1/2}\|\partial_x^5w^{(0)}\|_{L^\infty}^{1/2}
\leq\sqrt2 C_\chi N^5.
$$
Adding both bounds, $\|w^{(0)}\|_{C^{4,1/2}}\leq(5+\sqrt{2})C_\chi N^5\leq N^6$
as soon as $N$ is large enough.
 
We now turn to the Hilbert transform at the origin. Note that $w^{(0)}(x)=N^\varepsilon\varphi_0(N^{-\varepsilon}x)$, where $\varphi_0$ is the profile of Definition~\ref{defperfiles}. Since the Hilbert transform commutes with dilations, $Hw^{(0)}(x)=N^\varepsilon(H\varphi_0)(N^{-\varepsilon}x),$ so that $Hw^{(0)}(0)=N^\varepsilon H\varphi_0(0).$

As $N$ has been chosen large enough so that Lemma~\ref{estimaporabajoHilbert} applies with $M=N$, we have that $1/2\leq H\varphi_0(0)\leq3/2$. Hence
$$
\frac{N^\varepsilon}{2}\leq Hw^{(0)}(0)\leq\frac{3N^\varepsilon}{2}.
$$
The last assertion of the lemma is now immediate after taking $\alpha=1-9\varepsilon$.
\end{proof}

\begin{proof}[Proof of Theorem~\ref{construccionblowup}]
Fix $\delta\in(0,1)$ and choose $\varepsilon=\delta/9$. Pick $N$ big enough. To be precise, we need $N>1$ to be at least as large as in Theorem~\ref{thmIterativo} and Lemma~\ref{arrancar_ind}, and it also has to be large enough so that \begin{equation}\label{eq:crecimientoM}
        \exp\left(\frac{\varepsilon x^{\varepsilon/2}}{4}\right)\geq2x\qquad\text{for every }x\geq N,
    \end{equation}
    which is possible since $e^{\varepsilon x^{\varepsilon/2}/4}/x\to\infty$ as $x\to\infty$.

We claim that there exist a sequence of real numbers $\{M_j\}_{j\geq0}$ and two sequences $\{w^{(j)}\}_{j\geq0}$ and $\{F^{(j)}\}_{j\geq0}$ in $C^\infty(\R\times[0,1])$ such that, writing $t_j=1-M_j^{-\varepsilon/2}$, the following properties hold for every $j\geq0$:
\begin{enumerate}
    \item[(a)] $w^{(j)}$ solves Equation~\eqref{eq:DGF} with forcing term $F^{(j)}$, the function $w^{(j)}(\cdot,t)$ is odd and $\supp\bigl(w^{(j)}(\cdot,t)\bigr)\cup \supp\bigl(F^{(j)}(\cdot,t)\bigr)\subseteq[-1,1]$ for every $t\in[0,1]$. Moreover, $F^{(j)}\in C\bigl([0,1];C^{0,1-9\varepsilon}(\R)\bigr)$.
    \item[(b)] For every $t\in[t_j,1]$,
    $$
    \|w^{(j)}(\cdot,t)\|_{C^{4,1/2}}\leq M_j^6
    \qquad\text{and}\qquad
    \frac{M_j^\varepsilon}{4}\leq Hw^{(j)}(0,t)\leq2M_j^\varepsilon.
    $$
    \item[(c)] $M_j\geq2^jN$.
    \item[(d)] If $j\geq1$, $\bar w^{j}=w^{(j)}-w^{(j-1)}$ and $\bar F^{j}=F^{(j)}-F^{(j-1)}$ satisfy
    $\bar w^{j}(\cdot,t)=\bar F^{j}(\cdot,t)=0$ for every $t\in[0,t_{j-1}],$
    and
    $$
    \sup_{t\in[0,1]}\|\bar F^{j}(\cdot,t)\|_{C^{0,1-9\varepsilon}}\leq C_{\chi,\varepsilon}M_j^{-2\varepsilon}.
    $$
\end{enumerate}
We argue by induction on $j$. For $j=0$, set $M_0=N$ and let $w^{(0)}$ and $F^{(0)}$ be as in Lemma~\ref{arrancar_ind}, that is,
$$
w^{(0)}(x,t)=-M_0^\varepsilon\chi(M_0^{1/2-\varepsilon}x)\sin(M_0^{1-\varepsilon}x)
$$
and
$$
F^{(0)}(x,t)=\left[u(w^{(0)})\partial_xw^{(0)}-w^{(0)}Hw^{(0)}\right](x,t).
$$
Lemma~\ref{arrancar_ind} yields Items~(a) and~(b), while Item~(c) is trivial and Item~(d) is vacuous.
 
Assume now that $M_j$, $w^{(j)}$ and $F^{(j)}$ have been constructed and satisfy Items~(a),~(b) and~(c). Observe that $M_j$ is also big enough, and Items~(a) and~(b) are exactly the hypotheses of Theorem~\ref{thmIterativo} with $w^{(j)}$ as background and $F^{(j)}$ as its forcing term. We may therefore define
$$
M_{j+1}^{1/\varepsilon}=\exp\left(\int_{1-M_j^{-\varepsilon/2}}^1Hw^{(j)}(0,s)\dd s\right)
$$
and apply Theorem~\ref{thmIterativo}. This yields two functions $\bar w^{j+1}$ and $\bar F^{j+1}$ in $C^\infty(\R\times[0,1])$ such that $\bar w^{j+1}(\cdot,t)$ is odd for every $t\in[0,1]$, the forcing $\bar F^{j+1}$ belongs to $C\bigl([0,1];C^{0,1-9\varepsilon}(\R)\bigr)$, and $w^{(j+1)}=w^{(j)}+\bar w^{j+1}$
solves Equation~\eqref{eq:DGF} with forcing term $F^{(j+1)}=F^{(j)}+\bar F^{j+1}.$

Moreover, $\supp\bigl(w^{(j+1)}(\cdot,t)\bigr)\cup \supp\bigl(F^{(j+1)}(\cdot,t)\bigr)\subseteq[-1,1]$ for every $t\in[0,1]$, which completes the proof of Item~(a) for $j+1$.

Estimates~\eqref{eq:reg_pert} and~\eqref{HIPOTESISTHEN0} give that, for every $t\in[1-M_{j+1}^{-\varepsilon/2},1]=[t_{j+1},1]$,
$$
\|w^{(j+1)}(\cdot,t)\|_{C^{4,1/2}}\leq M_{j+1}^6
\qquad \text{and}\qquad
\frac{M_{j+1}^\varepsilon}{4}\leq Hw^{(j+1)}(0,t)\leq2M_{j+1}^\varepsilon,
$$
which is Item~(b) for $j+1$, while Property~\eqref{eq:antesdetiempo} and Inequality~\eqref{eq:estimafuerza} give $\bar w^{j+1}(\cdot,t)=\bar F^{j+1}(\cdot,t)=0$ for every $t\in[0,t_j]$ and
$$
\sup_{t\in[0,1]}\|\bar F^{j+1}(\cdot,t)\|_{C^{0,1-9\varepsilon}}\leq C_{\chi,\varepsilon}M_{j+1}^{-2\varepsilon}
$$
respectively, which is Item~(d) for $j+1$. 

Finally, the lower bound in Item~(b) for $w^{(j)}$ implies
$$
M_{j+1}^{1/\varepsilon}
=\exp\left(\int_{1-M_j^{-\varepsilon/2}}^1Hw^{(j)}(0,s)\dd s\right)
\geq\exp\left(\frac{M_j^{\varepsilon/2}}{4}\right)
$$
and hence, by Inequality~\eqref{eq:crecimientoM},
$$
M_{j+1}\geq\exp\left(\frac{\varepsilon M_j^{\varepsilon/2}}{4}\right)\geq2M_j\geq2^{j+1}N.
$$
This is Item~(c) for $j+1$, which closes the induction.
 
Let us make the following observations. First, as $M_{j+1}\geq 2M_j$, the sequence $\{t_j\}_{j\geq0}$ is increasing in $(0,1)$, with $t_j\nearrow1$. Second, if $i\geq j$, then $w^{(i)}-w^{(j)}=\sum_{l=j+1}^{i}\bar w^{l}$ and, by Item~(d), each $\bar w^l$ with $l\geq j+1$ vanishes for $t\in[0,t_{l-1}]\supseteq[0,t_j]$. The same applies to the forcing terms, and therefore
\begin{equation}\label{eq:estabilizacion}
    w^{(i)}(\cdot,t)=w^{(j)}(\cdot,t),
    \quad
    F^{(i)}(\cdot,t)=F^{(j)}(\cdot,t)
    \qquad\text{for every }t\in[0,t_j]\text{ and every }i\geq j.
\end{equation}
Third, observe that by Item~(a):
\begin{equation}\label{eq:soportefuerza}
    \supp\bigl(F^{(j)}(\cdot,t)\bigr)\subseteq[-1,1]
    \qquad\text{for every }t\in[0,1]\text{ and every }j\geq0.
\end{equation}
For every $x\in\R$ and $t\in[0,1)$ we define
$$
w(x,t)=\lim_{j\to\infty}w^{(j)}(x,t).
$$
This limit exists, as the sequence $\{w^{(i)}(x,t)\}_{i\geq 0}$ is eventually constant. Indeed, since $t_j\nearrow1$, for every $t\in[0,1)$ there exists $j$ with $t\leq t_j$, and by the first equation in Property~\eqref{eq:estabilizacion}, the sequence $\{w^{(i)}(x,t)\}_{i\geq0}$ is constant for $i\geq j$. In particular,
\begin{equation}\label{eq:wigualwj}
    w=w^{(j)}\qquad\text{in }\R\times[0,t_j],\ \text{for every }j\geq0.
\end{equation}
With regard to the forcing term, note that by Items~(a),~(c) and~(d), every term of the series $F^{(0)}+\sum_{j\geq1}\bar F^{j}$ belongs to the Banach space $C\bigl([0,1];C^{0,1-9\varepsilon}(\R)\bigr)$, and
\begin{equation}\label{eq:seriefuerza}
    \sum_{j=1}^\infty\sup_{t\in[0,1]}\|\bar F^{j}(\cdot,t)\|_{C^{0,1-9\varepsilon}}
    \leq C_{\chi,\varepsilon}\sum_{j=1}^\infty M_j^{-2\varepsilon}
    \leq\frac{C_{\chi,\varepsilon}}{N^{2\varepsilon}}\sum_{j=1}^\infty2^{-2\varepsilon j}
    =\frac{C_{\chi,\varepsilon}}{N^{2\varepsilon}(2^{2\varepsilon}-1)}<\infty.
\end{equation}
Hence the series converges absolutely in $C\bigl([0,1];C^{0,1-9\varepsilon}(\R)\bigr)$, and we can define
$$
F=F^{(0)}+\sum_{j=1}^\infty\bar F^{j}=\lim_{j\to\infty}F^{(j)}
\quad\text{in }C\bigl([0,1];C^{0,1-9\varepsilon}(\R)\bigr).
$$
Note that $F$ is defined for every $t\in[0,1]$, including the blow-up time. Since convergence in this space implies uniform convergence in $\R\times[0,1]$, we also have $F(x,t)=\lim_{j\to\infty}F^{(j)}(x,t)$ for every $(x,t)\in\R\times[0,1]$. Together with the second equation of Property~\eqref{eq:estabilizacion} and Property~\eqref{eq:soportefuerza}, this gives
\begin{equation}\label{eq:FigualFj}
    F=F^{(j)}\qquad\text{in }\R\times[0,t_j],\ \text{for every }j\geq0,
\end{equation}
and $\supp\bigl(F(\cdot,t)\bigr)\subseteq[-1,1]$ for every $t\in[0,1]$.
 
Observe that $w$ is a classical solution with forcing $F$ on $[0,1)$. Indeed, let $c\in(0,1)$ and choose $j$ such that $t_j\geq c$. By Equations~\eqref{eq:wigualwj} and~\eqref{eq:FigualFj}, $w=w^{(j)}$ and $F=F^{(j)}$ in $\R\times[0,c]$. Since $w^{(j)},F^{(j)}\in C^\infty(\R\times[0,1])$ and $c\in(0,1)$ is arbitrary, we deduce that $w,F\in C^\infty(\R\times[0,1)).$

Moreover, as $w^{(j)}$ solves Equation~\eqref{eq:DGF} with forcing term $F^{(j)}$, it follows that the restriction of $w$ to $\R\times[0,c]$ is a classical solution on $[0,c]$ with forcing $F$. Finally, by Item~(a), Equation~\eqref{eq:wigualwj} and the previous step,
$$
\supp\bigl(w(\cdot,t)\bigr)\cup\supp\bigl(F(\cdot,t)\bigr)\subseteq[-1,1]
\qquad\text{for every }t\in[0,1).
$$
Therefore, $w$ is a classical solution to Equation~\eqref{eq:DGF} with forcing $F$ on $[0,1)$, with $K=[-1,1]$.
 
Regarding the regularity of the force, we have that $F\in C\bigl([0,1];C^{0,1-\delta}(\R)\bigr)$. Moreover, by Inequality~\eqref{eq:seriefuerza},
$$
\sup_{t\in[0,1]}\|F(\cdot,t)\|_{C^{0,1-9\varepsilon}}
\leq\|F^{(0)}\|_{C^{0,1-9\varepsilon}}
+\frac{C_{\chi,\varepsilon}}{N^{2\varepsilon}(2^{2\varepsilon}-1)}
<\infty.
$$ 
Finally, we show that the solution blows up. By Equation~\eqref{eq:wigualwj}, $w(\cdot,t_j)=w^{(j)}(\cdot,t_j)$, and by Item~(b),
$$
\|Hw(\cdot,t_j)\|_{L^\infty}
\geq Hw(0,t_j)
=Hw^{(j)}(0,t_j)
\geq\frac{M_j^\varepsilon}{4},
$$
while
$$
\|w(\cdot,t_j)\|_{C^1}
\leq\|w^{(j)}(\cdot,t_j)\|_{C^{4,1/2}}
\leq M_j^6.
$$
In particular, $w(\cdot,t_j)$ is not identically zero. Applying Lemma~\ref{estimalogH} with $R=1$ we obtain
\begin{equation}\label{logblowupfinal}
\|Hw(\cdot,t_j)\|_{L^\infty}
\leq C\|w(\cdot,t_j)\|_{L^\infty}
\left[1+\log\left(\frac{\|w(\cdot,t_j)\|_{C^1}}{\|w(\cdot,t_j)\|_{L^\infty}}\right)\right]
\end{equation}
for every $j\in\mathbb{N}$ and some fixed $C>0$.

We now claim that $\|w(\cdot,t_j)\|_{L^\infty}\geq1$ for infinitely many $j\in\N$. Otherwise, there would exist $j_0\in\mathbb{N}$ sufficiently large such that $\|w(\cdot,t_j)\|_{L^\infty}<1$ for every $j\geq j_0$. For each $j\in\mathbb{N}$, consider the function
$$
y_j(x)=Cx\left[1+\log\left(\frac{M_j^6}{x}\right)\right].
$$
Since $y_j$ is increasing on $(0,1]$, Inequality~\eqref{logblowupfinal} would imply $M_j^\varepsilon\leq 4C(1+\log(M_j^6))$
for every $j\geq j_0$, which is impossible since $M_j\to\infty$ as $j\to\infty$.

Thus, there exists a subsequence $\{t_{n_j}\}_{j\in\mathbb{N}}$ of $\{t_j\}_{j\in\mathbb{N}}$ such that, by~\eqref{logblowupfinal},
$$
\|w(\cdot,t_{n_j})\|_{L^\infty}
\geq \frac{cM_{n_j}^\varepsilon}{1+\log(M_{n_j})}
$$
for every $j\in\mathbb{N}$, where $c>0$ is constant. Since $M_j\to\infty$ as $j\to\infty$, we conclude that
$$
\limsup_{t\nearrow1}\|w(\cdot,t)\|_{L^\infty}=\infty,
$$
which completes the proof. 
\end{proof}

\section*{Acknowledgments} 
This work was supported in part by the Spanish Ministry of Science and Innovation through the ``Severo Ochoa Programme for Centres of Excellence in R\&D (CEX2023-001347-S)'' and grant PID2023-152878NB-I00. J.A.F.M. acknowledges partial support received from a JAE Intro ICU 2025 fellowship awarded by the Spanish National Research Council (CSIC).

\section*{AI disclosure} 
Claude Opus 5.5 was used in the preparation of this manuscript for the following purposes: generating the TikZ code for Figure 1, preparing the BibTeX entries for the bibliography and translating several sentences from a preliminary draft originally written in Spanish.

\bibliographystyle{plainnat}
\bibliography{bibliografiap}

\end{document}